%% file: main.tex
\documentclass[11pt,reqno]{amsart}
\usepackage[colorlinks=true,
pdfstartview=FitV, linkcolor=cyan, citecolor=magenta,
urlcolor=blue]{hyperref}
\usepackage{amsmath,amsfonts,latexsym,amssymb}
\usepackage{mathrsfs}
\usepackage[latin1]{inputenc}
\usepackage[T1]{fontenc}
\usepackage{ae,aecompl}
\usepackage{braket}
\usepackage{comment}
\usepackage{color}
\usepackage{graphicx}
\usepackage{subfig}
\usepackage{bbm}

\usepackage{enumitem}

\newtheorem{theorem}{Theorem}[section]
\newtheorem{lemma}[theorem]{Lemma}
\newtheorem{proposition}[theorem]{Proposition}
\newtheorem{corollary}[theorem]{Corollary}

\newtheorem*{thm*}{\protect\theoremname}
\theoremstyle{definition}

\newtheorem{remark}[theorem]{Remark}

\newtheorem*{cor*}{\protect\corollaryname}

\long\def\@savemarbox#1#2{\global\setbox#1\vtop{\hsize\marginparwidth 
  \@parboxrestore\tiny\raggedright #2}}
\newcommand{\Cc}{{\mathcal C}}

\newcommand{\PSL}{\mathsf{PSL}}
\newcommand{\PGL}{\mathsf{PGL}}

\newcommand{\abs}[1]{\left|#1\right|}

\newcommand{\Oc}{\mathcal O}

\DeclareMathOperator{\Aut}{\mathsf{Aut}}
\DeclareMathOperator{\id}{\mathrm{id}}

\DeclareMathOperator{\Cb}{\mathbb{C}}

\DeclareMathOperator{\Hb}{\mathbb{H}}

\DeclareMathOperator{\Kb}{\mathbb{K}}

\DeclareMathOperator{\Nb}{\mathbb{N}}
\DeclareMathOperator{\Pb}{\mathbb{P}}
\DeclareMathOperator{\Rb}{\mathbb{R}}

\DeclareMathOperator{\Usf}{\mathsf{U}}
\DeclareMathOperator{\Psf}{\mathsf{P}}

\DeclareMathOperator{\Gsf}{\mathsf{G}}

\DeclareMathOperator{\Fc}{\mathcal{F}}

\DeclareMathOperator{\SL}{\mathsf{SL}}
\newcommand{\norm}[1]{\left\|#1\right\|}

\DeclareMathOperator{\Span}{Span}

\DeclareMathOperator{\Gr}{\mathrm{Gr}}

\usepackage{ mathrsfs }

\providecommand{\corollaryname}{Corollary}
\providecommand{\theoremname}{Theorem}

\begin{document}

\title{A rigidity theorem for complex Kleinian groups}
\author[Canary]{Richard Canary}
\address{University of Michigan}
\author[Zhang]{Tengren Zhang}
\address{National University of Singapore}
\author[Zimmer]{Andrew Zimmer}
\address{University of Wisconsin-Madison}
\thanks{Canary was partially supported by grant  DMS-2304636  from the National Science Foundation.
Zhang was partially supported by the NUS-MOE grants A-8001950-00-00 and A-8000458-00-00. Zimmer was partially supported by a Sloan research fellowship and grants DMS-2105580 and DMS-2452068 from the
National Science Foundation.}

\begin{abstract} 
Farre, Pozzetti and Viaggi proved that any $(d-k)$-hyperconvex subgroup of $\mathsf{PSL}(d,\mathbb C)$  is virtually isomorphic to a convex cocompact Kleinian group
and that its $\alpha_k$-critical exponent is at most 2. We show that a $(d-k)$-hyperconvex subgroup is isomorphic to a uniform lattice in $\mathsf{PSL}(2,\mathbb C)$
if and only if its $\alpha_k$-critical exponent is exactly  2.  Furthermore, we show that if a  strongly irreducible $(d-k)$-hyperconvex subgroup  has
$\alpha_k$-critical exponent $2$, then it is the image of a uniform lattice in $\mathsf{PSL}(2,\mathbb C)$ by an irreducible representation of 
$\mathsf{PSL}(2,\mathbb C)$  into $\mathsf{PSL}(d,\mathbb C)$.
\end{abstract}

\maketitle

\begin{center}
{\em Dedicated to Pepe Seade on the occasion of his seventieth birthday}
\end{center}

\tableofcontents

\input{intro}

\input{background}

\input{ahlfors}


\input{lattice}

\input{appendix}

\input{bib}

\end{document}

%% file: intro.tex

\section{Introduction}
In this paper, we continue the study of $(d-k)$-hyperconvex subgroups of $\mathsf{PSL}(d,\mathbb C)$ initiated by
Pozzetti--Sambarino--Wienhard \cite{PSW1} and Farre--Pozzetti--Viaggi \cite{FPV}. A subgroup of $\mathsf{PSL}(d,\mathbb C)$ is
$(d-k)$-hyperconvex for some $k\in\{1,\ldots, d-1\}$ if it is Anosov with respect to the simple roots $\alpha_{k-1}$, $\alpha_{k}$ and $\alpha_{k+1}$
and triples of points in its limit set satisfy a specific transversality condition (see Section~\ref{sec: hyperconvex} for the definition). The most basic examples of $k$-hyperconvex groups are obtained by considering the image of
a convex cocompact subgroup of $\PSL(2,\mathbb C)$ by an irreducible representation of $\PSL(2,\mathbb C)$
into $\mathsf{PSL}(d,\mathbb C)$.

Farre--Pozzetti--Viaggi showed that any $(d-k)$-hyperconvex subgroup of $\mathsf{PSL}(d,\mathbb C)$ is virtually isomorphic
to a convex cocompact subgroup of $\mathsf{PSL}(2,\mathbb C)$. They further showed that the
$\alpha_k$-critical exponent $\delta^{\alpha_k}(\Gamma)$ of a $k$-hyperconvex subgroup $\Gamma$  of $\mathsf{PSL}(d,\mathbb C)$ is at most 2. 
Recall that the $\alpha_k$-critical exponent $\delta^{\alpha_k}(\Gamma)$ is the exponential growth rate of the number of elements $\gamma\in\Gamma$
so that 
$$\alpha_k(\kappa(\gamma))=\log \frac{\sigma_k(\gamma)}{\sigma_{k+1}(\gamma)}\le T,$$
where $\sigma_1(\gamma) \ge \cdots \ge \sigma_d(\gamma)$ are the singular values of $\gamma$. 

It then follows from work of Pozzetti--Sambarino--Wienhard \cite{PSW1}, that if $\Gamma$ is $(d-k)$-hyperconvex, then the Hausdorff dimension of
the limit set $\Lambda_k(\Gamma)$ of the action of $\Gamma$ on the Grassmanian $\mathrm{Gr}_k(\mathbb C^d)$ is exactly
$\delta^{\alpha_k}(\Gamma)\le 2$. If $\Gamma$ is virtually isomorphic to a lattice in $\mathsf{PSL}(2,\mathbb C)$, then $\Lambda_k(\Gamma)$ is 
homeomorphic to a two-sphere, so we must have $\delta^{\alpha_k}(\Gamma)=2$.  In this paper, we prove that the converse also holds.

\begin{theorem}
\label{first result}
If $\Gamma\subset \PSL(d,\mathbb C)$ is $(d-k)$-hyperconvex, then 
$\delta^{\alpha_{k}}(\Gamma)=2$ if and only if  $\Gamma$ is virtually isomorphic to a uniform lattice in $\PSL(2,\mathbb C)$.
\end{theorem}

One may view this as a generalization of the fact, due to Sullivan \cite{sullivan} and Tukia \cite{tukia}, that  the critical exponent of a convex cocompact 
subgroup of $\PSL(2,\mathbb C)$ has Hausdorff dimension two if and only if it is a uniform lattice.

\medskip

The only known examples where $\Gamma$ is $(d-k)$-hyperconvex, strongly irreducible and virtually isomorphic to a lattice in $\PSL(2,\mathbb C)$ are obtained from lattices
in $\PSL(2,\mathbb C)$ by applying an irreducible representation of $\PSL(2,\mathbb C)$ into $\PSL(d,\mathbb C)$. We show that these are the only examples.

\begin{theorem}
\label{rigidity result}
If $\Gamma\subset \PSL(d,\mathbb C)$ is $(d-k)$-hyperconvex, strongly irreducible and 
$\delta^{\alpha_{k}}(\Gamma)=2$, then 
there exist a uniform lattice $\Gamma' \subset \PSL(2,\Cb)$ and  an irreducible representation $\tau:\PSL(2,\mathbb C)\to\PSL(d,\mathbb C)$ so
that  $\Gamma=\tau(\Gamma')$.
\end{theorem}

Since the property of being $(d-k)$-hyperconvex is stable (see \cite[Prop.\ 6.2]{PSW1}), one may view our result as a global analogue of the local rigidity result of  
Raghunathan \cite{raghunathan} (see also Menal--Ferrer and Porti \cite{MP,porti}) who proved that if $\Gamma$ is  a uniform lattice in $\mathsf{PSL}(2,\mathbb C)$, then 
$\tau|_\Gamma$ is locally rigid in the character variety ${\bf X}(\Gamma,\mathsf{PSL}(d,\mathbb C))$
if $\tau:\mathsf{PSL}(2,\mathbb C)\to\mathsf{PSL}(d,\mathbb C)$ is an  irreducible representation.

As a corollary to Theorem~\ref{rigidity result} we have the following.

\begin{corollary}
\label{mostow-like}
If $\Gamma' \subset \PSL(2,\Cb)$ is a uniform lattice and $\rho : \Gamma' \rightarrow \PSL(d,\Cb)$ is a strongly irreducible representation with finite kernel and $(d-k)$-hyperconvex image, then $\rho$ extends to an irreducible representation $\PSL(2,\Cb) \rightarrow \PSL(d,\Cb)$. 
\end{corollary} 

Any uniform lattice in $\PSL(2,\Cb)$ is strongly irreducible and 1-hyperconvex, hence Corollary \ref{mostow-like} can be viewed as a generalization of Mostow's 
rigidity theorem \cite{mostow}  for hyperbolic 3-manifolds. 

\medskip

We also note that without irreducibility, there is considerable less rigidity. For instance, fix a uniform lattice $\Gamma_0 \subset \mathsf{SL}(2,\Cb)$ and \textbf{any} representation $\rho : \Gamma_0 \rightarrow \mathsf{SL}(d_0, \Cb)$. Next let $\tau_d : \mathsf{SL}(2,\Cb) \rightarrow \mathsf{SL}(d,\Cb)$ denote the standard irreducible representation and define $\rho_d : \Gamma_0 \rightarrow \mathsf{PSL}(d+d_0, \Cb)$ by 
$$
\rho_d(\gamma) = \begin{bmatrix} \rho(\gamma) & \\ & \tau_d(\gamma) \end{bmatrix}.
$$
Then for any $k \in \Nb$ and $d$ sufficiently large (depending on $\rho$ and $k$) the image $\Gamma : = \rho_d(\Gamma_0)$ is $(d-k)$-hyperconvex and has $\delta^{\alpha_{k}}(\Gamma)=2$.

\subsection*{The context} We would like to place our results into two natural contexts. One may naturally ask what restrictions on an Anosov group place restrictions on its
isomorphism type. One early example of a result of this form is due to Guichard \cite{guichard} who showed that if $\Gamma$ is a $\Psf_1$-Anosov subgroup of 
$\mathsf{PSL}(d,\mathbb R)$ which is isomorphic to a surface group, so that any $d$ lines in $\Lambda_1(\Gamma)$ span $\mathbb R^d$, then $\Gamma$
is the image of a Hitchin representation.
To our knowledge, the papers of Pozzetti--Sambarino--Wienhard \cite{PSW1} and Canary--Tsouvalas \cite{canary-tsouvalas} are the first general investigations of this 
question. In the $\mathsf{PSL}(d,\mathbb R)$ setting, Sambarino asked whether every Borel Anosov subgroup of $\mathsf{PSL}(d,\mathbb R)$ is virtually isomorphic
to a free group or a surface group. This surprising question is now known to have a positive answer whenever
$d=2,\ 3,\ 4,\ 5, \text{ or }6\ \mathrm{mod}\ 8$, due to work of Dey  \cite{dey-sambarino} and Tsouvalas \cite{tsouvalas-sambarino}.

Another natural context is the theory of complex Kleinian groups which was developed by Jos\'e Seade and his collaborators. A complex Kleinian group
is a discrete subgroup of $\mathsf{PSL}(d,\mathbb C)$. The theory of complex Kleinian groups developed independently and contemporaneously with
the theory of Anosov representations. It seems to us that the stage is now set for these two fields to interact profitably. We recommend that book
 \cite{CNS} as an introduction to the field of complex Kleinian groups.

\subsection*{Outline of proofs} The key tool in the proof of Theorem \ref{first result} is a deeper analysis of the limit sets of hyperconvex groups.
We use the technology developed in \cite{CZZ2} and \cite{CZZ3} to show that if $\Gamma\subset\mathsf{PSL}(d,\mathbb C)$ is $(d-k)$-hyperconvex,
then the $\alpha_k$-Patterson-Sullivan measure $\mu$ on $\Lambda_k(\Gamma)$  is $\delta:=\delta^{\alpha_{k}}(\Gamma)$-regular, 
i.e. there exist  $D\ge 1$ and $t_0>0$ such that for every $x\in \Lambda_k(\Gamma)$ and every $0<t<t_0$, we have
\[\frac{1}{D}t^\delta\le \mu (B(x,t))\le Dt^\delta,\]
where $B(x,t)$ is the open ball of radius $t$ about $x$ (see Theorem \ref{prop: Ahlfors regularity}). 
Dey--Kim--Oh \cite{DKO}, see also Dey--Kapovich\cite{DK}, previously showed that if  $\Gamma$ is $\Psf_\theta$-Anosov, $\phi\in \mathfrak a_\theta^*$ is symmetric and 
$\delta^\phi(\Gamma)<\infty$, then the $\phi$-Patterson-Sullivan measure on $\Lambda_\theta(\Gamma)$ is $\delta^\phi(\Gamma)$-Ahlfors regular with respect
to a Finsler premetric on $\mathcal F_\theta$. In our case, $\alpha_k$ is not symmetric and the Ahlfors regularity is with respect to the Riemannian metric
on $\mathrm{Gr}_k(\mathbb C^d)$.

The following result is a simple consequence of Ahlfors regularity.

\begin{theorem}
\label{finite Hausdorff measure}
If $\Gamma\subset \PSL(d,\mathbb C)$ is $(d-k)$-hyperconvex, then $\Lambda_k(\Gamma)$ 
has finite, non-zero $\delta^{\alpha_k}(\Gamma)$-dimensional Hausdorff measure.
\end{theorem}

Farre--Pozzetti--Viaggi  \cite[Thm.\ B]{FPV}  showed that if $\Gamma\subset\mathsf{PSL}(d,\mathbb C)$  is $(d-k)$-hyperconvex and $\Gamma$ 
is not virtually isomorphic to a uniform lattice in 
$\mathsf{PSL}(2,\mathbb C)$, then the limit set $\Lambda_k(\Gamma)$ of $\Gamma$ has 2-dimensional Hausdorff measure zero.
Theorem \ref{first result} follows immediately from their result and Theorem \ref{finite Hausdorff measure}.

We note that in the proof of Theorem \ref{finite Hausdorff measure} one may use  standard properties
of the exterior power representation to reduce the general case to the case where $k=1$.

\medskip

One key input in the proof of Theorem \ref{rigidity result} is the fact, which follows from work of Pozzetti--Sambarino--Wienhard \cite{PSW1} 
(see also Zhang--Zimmer \cite{ZZ}), that if $\Gamma$ is
$(d-k)$-hyperconvex and $\partial\Gamma$ is a two-sphere, then  the limit set $\Lambda_{k-1,k,k+1}$ of $\Gamma$ in the partial
flag variety $\mathcal F_{k-1,k,k+1}$ is  a $\mathcal C^1$-smooth submanifold. 
(Elements of $\mathcal F_{k-1,k,k+1}$ are partial flags with subspaces in dimensions $k-1$, $k$ and $k+1$).
We then observe that  $\Lambda_{k-1,k,k+1}(\Gamma)$ is a complex submanifold of $\mathrm{Gr}_k(\mathbb C^d)$, and hence, 
by Chow's theorem \cite{chow}, a complex algebraic 
subvariety of  $\mathcal F_{k-1,k,k+1}$. This implies that the real Zariski closure $\mathsf G$ of $\Gamma$ acts by biholomorphisms on 
$\Lambda_{k-1,k,k+1}(\Gamma)$. Then, since $\Lambda_{k-1,k,k+1}(\Gamma)$ is a two-sphere and hence biholomorphic to $\mathbb{CP}^1$, we obtain a surjective homomorphism from  the identify component $\mathsf G^0$ of $\mathsf G$ to $\mathsf{PSL}(2,\mathbb C)$, 
the biholomorphism group of $\mathbb{CP}^1$.  We use Lie theoretic arguments to prove that $\mathsf G$ is simple and then  geometric
arguments to conclude that $\mathsf G$ is isomorphic to $\mathsf{PSL}(2,\mathbb C)$.

%% file: background.tex

\section{Background}

\subsection{The Jordan and Cartan projections}\label{The Jordan and Cartan projections}
Let $\Kb$ be either the field of real or complex numbers, and let $(e_1,\dots,e_d)$ be the standard basis of $\Kb^d$. Recall that the standard {\em Cartan subspace} for
$\PSL(d,\Kb)$, denoted $\mathfrak a$, is the real vector space of real-valued, $d\times d$, diagonal, traceless matrices:
$$\mathfrak a=\{\mathrm{diag}(A_1,\ldots, A_d)\in\mathfrak{sl}(d,\Rb)\ |\ A_1+\cdots +A_d=0\}.$$
For each $k\in\Delta:=\{1,\dots,d-1\}$, the \emph{$k$-th restricted simple root} is the linear map
\[\alpha_k:\mathfrak a\to\Rb\] 
given by $\alpha_k(A)=A_k-A_{k+1}$. The standard {\em positive Weyl chamber} is the subset where all the
simple roots are non-negative:
$$\mathfrak a^+=\{A\in \mathfrak a\ |\ \alpha_k(A)\ge 0\text{ for all }k\in\Delta\}.$$

Given $g\in\PSL(d,\Kb)$, let  
\[\lambda_1(g)\ge\dots\ge\lambda_d(g)\] 
denote the modulus of the generalized eigenvalues of some lift $\tilde g$ of $g$ in $\SL(d,\Kb)$, and let 
\[\sigma_1(g)\ge\dots\ge\sigma_d(g)\]
denote the singular values of $\tilde g$ (with respect to the standard Hermitian inner product on $\Kb^d$). The {\em Jordan projection} and \emph{Cartan projection} 
\[\nu,\kappa:\PSL(d,\mathbb K)\to \mathfrak a^+\] 
are respectively given by
\[\nu(g)=\mathrm{diag}(\log\lambda_1(g),\ldots,\log\lambda_d(g))\quad\text{and}\quad\kappa(g)=\mathrm{diag}(\log\sigma_1(g),\ldots,\log\sigma_d(g)).\]
Recall that (see \cite[Cor.\ 6.34]{benoist-quint})
$$\lambda(g) = \lim \frac{1}{n} \kappa(g^n).$$

By the singular value decomposition theorem, any $g\in\PSL(d,\mathbb K)$ can be decomposed as $g=m_1\exp(\kappa(g))m_2$ for some $m_1,m_2 \in\mathsf{PSU}(d,\mathbb K)$. Furthermore, if $k\in\{1,\dots,d-1\}$ satisfies $\alpha_k(\kappa(g))>0$, then 
\[U_k(g):=m_1(\Span_{\Kb}( e_1,\ldots,e_k)) \in \Gr_k(\Kb^d)\] 
is well-defined (even though $m_1$ and $m_2$ are not uniquely specified by $g$). In this case, one can verify that $U_{d-k}(g^{-1})$ is also well-defined, and can be described as 
\[U_{d-k}(g^{-1})=m_2^{-1}(\Span_{\Kb}( e_{k+1},\ldots,e_d)).\]

\subsection{Flag manifolds and the Iwasawa decomposition}
Given a non-empty subset 
\[\theta=\{k_1,k_2,\dots,k_r\}\subset\Delta,\] 
the \emph{$\theta$-flag manifold}, denoted $\Fc_\theta=\Fc_\theta(\Kb^d)$, is the set of nested subspaces 
$x^{k_1}\subset\dots\subset x^{k_r}\subset\Kb^d$ so that $x^{k_j}$ has dimension $k_j$ for each $j$.
We often denote $x:=(x^{k_1},\dots,x^{k_r})$. Let $\Psf_\theta\subset\PSL(d,\Kb)$ be the stabilizer of the flag $x_0\in\Fc_\theta$ defined by $x_0^k=\Span(e_1,\dots,e_k)$ 
for all $k\in\theta$, where $e_1,\dots, e_d$ is the standard basis of $\Kb^d$. 
Since $\PSL(d,\Kb)$ acts transitively on the space of flags, we have the identification $\PSL(d,\Kb)/\mathsf P_\theta\cong\Fc_\theta$ as $\PSL(d,\Kb)$-spaces.

Let $\iota:\Delta\to\Delta$ be the map given by $\iota(k)=d-k$. We say that flags $x\in\Fc_\theta$ and $y\in \Fc_{\iota(\theta)}$ are \emph{transverse} if $x^k+y^{d-k}=\Kb^d$ for all $k\in\theta$. For every $x\in\Fc_\theta$, let $\mathcal Z_x\subset\Fc_{\iota(\theta)}$ denote the set of flags that are not transverse to $x$.

If $g\in\PSL(d,\Kb)$ satisfies $\alpha_k(\kappa(g))>0$ for all $k\in\theta$, then we may define $U_\theta(g)\in\Fc_\theta$ by $U_\theta(g)^k:=U_k(g)$ for all $k\in\theta$. The following is standard (see for example, \cite[Prop.\ 2.3]{CZZ3}).

\begin{proposition}\label{prop:characterizing convergence in general} Let $\theta\subset\Delta$ be non-empty. Suppose $x^+\in\Fc_\theta$, $x^- \in \Fc_{\iota(\theta)}$ and $(g_n)_{n=1}^\infty$ is a sequence in $\mathsf{PSL}(d,\Kb)$. Then the following are equivalent:  
\begin{enumerate}
\item $\alpha_k(\kappa(g_n)) \to \infty$ for every $k \in \theta$, $U_\theta(g_n) \rightarrow x^+$ and $U_{\iota(\theta)}(g_n^{-1}) \rightarrow x^-$.
\item $g_n(x) \to x^+$ for all $x \in \Fc_\theta \setminus \mathcal{Z}_{x^-}$, and this convergence is uniform on compact subsets of $\Fc_\theta \setminus  \mathcal{Z}_{x^-}$. 
\item $g_n^{-1}(x) \to x^-$ for all $x \in \Fc_{\iota(\theta)} \setminus  \mathcal{Z}_{x^+}$, and this convergence is uniform on compact subsets of $\Fc_{\iota(\theta)} \setminus \ \mathcal{Z}_{x^+}$. 
\item There are open sets $\mathcal{U}^+\subset\Fc_\theta$ and $\mathcal U^-\subset\Fc_{\iota(\theta)}$ such that $g_n(x) \to x^+$ for all $x \in \mathcal U^+$ and $g_n^{-1}(x) \to x^-$ for all $x \in \mathcal U^-$.
\end{enumerate}
\end{proposition}

By the Iwasawa decomposition theorem, we may decompose every $g\in\PSL(d,\Kb)$ uniquely as $g=me^An$ for some $m\in\mathsf{PSU}(d,\Kb)$, $A\in\mathfrak a$, and $n$ a unipotent upper triangular matrix in $\PSL(d,\Kb)$. Since $\mathsf{PSU}(d,\Kb)$ acts transitively on $\Fc_\Delta$, we may use this to define the \emph{Iwasawa cocycle} 
\[B:\PSL(d,\Kb)\times\Fc_{\Delta}\to\mathfrak a\]
as follows. Given any $g\in\PSL(d,\Kb)$ and any $x\in\Fc_\Delta$, choose $k\in\mathsf{PSU}(d,\Kb)$ such that $x=k \Psf_\Delta$. Then decompose $gk=me^An$. It is straightforward to verify that $A$ does not depend on the choice of $k$, so we may define $B(g,x):=A$. 

\subsection{Divergent, transverse and Anosov subgroups}\label{Divergent, transverse and Anosov subgroups}
For any non-empty $\theta\subset\Delta$, a subgroup $\Gamma\subset\PSL(d,\Kb)$ is 
\emph{$\mathsf{P}_\theta$-divergent} if for any $k\in\theta$ and any infinite sequence $\{\gamma_i\}$ in $\Gamma$, we have
\[\lim_{i\to\infty}\alpha_k(\kappa(\gamma_i))=\infty.\]

For any $\mathsf{P}_\theta$-divergent subgroup $\Gamma\subset\PSL(d,\Kb)$, we may define its \emph{$\theta$-limit set}
\[\Lambda_\theta(\Gamma):=\left\{\lim_{i\to\infty}U_\theta(\gamma_i)\in\Fc_\theta:(\gamma_i)_{i=1}^\infty\text{ is an infinite sequence in }
\Gamma\right\}.\] 
One can verify that $\Lambda_\theta(\Gamma)$ is a $\Gamma$-invariant subset of $\Fc_\theta$, see \cite[Prop.\ 3.1]{CZZ2}. If $\theta$ is symmetric, i.e. $\theta=\iota(\theta)$, a $\mathsf{P}_\theta$-divergent subgroup $\Gamma\subset\PSL(d,\Kb)$ is called \emph{$\mathsf{P}_\theta$-transverse} if 
$\Lambda_\theta(\Gamma)\subset\Fc_\theta$ is a transverse subset, i.e. any pair of distinct flags in $\Lambda_\theta(\Gamma)$ are transverse.

Recall that the action of a group $\Gamma$ by homeomorphisms on a compact metrizable topological space $M$ is called a \emph{convergence group action} if for any infinite sequence $(\gamma_i)_{i=1}^\infty$ in $\Gamma$, there is a subsequence $(\gamma_{i_k})_{k=1}^\infty$ and a pair of points $a,b\in M$ such that $\gamma_{i_k}(x)\to a$ for all $x\in M-\{b\}$. If $\Gamma$ acts on $M$ as a convergence group, a point $x\in M$ is a \emph{conical point} if there exists a sequence $(\gamma_i)_{i=1}^\infty$ in $\Gamma$ and a distinct pair of points $a,b\in M$ such that $\gamma_i(x)\to a$ and $\gamma_i(y)\to b$ for all $y\in M-\{x\}$. 

Using Proposition \ref{prop:characterizing convergence in general}, one can verify that if $\Gamma\subset\PSL(d,\Kb)$ is $\mathsf{P}_\theta$-transverse, then it acts on $\Lambda_\theta(\Gamma)$ as a convergence group, see \cite[Prop.\ 3.3]{CZZ2} or \cite[Sec.\ 5.2]{KLP}. We say that a $\mathsf{P}_\theta$-transverse subgroup $\Gamma\subset\PSL(d,\Kb)$ is \emph{$\mathsf{P}_\theta$-Anosov} if every point in $\Lambda_\theta(\Gamma)$ is conical. Equivalently, a $\Psf_\theta$-transverse subgroup $\Gamma\subset\PSL(d,\Kb)$ is $\mathsf{P}_\theta$-Anosov if and only if $\Gamma$ is word hyperbolic and there exists a $\Gamma$-equivariant homeomorphism between $\Lambda_\theta(\Gamma)$ and the Gromov boundary $\partial\Gamma$ of $\Gamma$.

To avoid cumbersome notation, we often say a group is $\Psf_k$-Anosov instead of $\Psf_{\{k,d-k\}}$-Anosov. 

The following is a characterization of Anosov representations due to Kapovich--Leeb--Porti \cite[Thm.\ 1.5]{KLP-morse} and Bochi--Potrie--Sambarino \cite[Prop.\ 4.6 and 4.9]{BPS}.

\begin{theorem}\label{Anosov singular values}
A subgroup $\Gamma\subset\PSL(d,\Kb)$ is $\Psf_k$-Anosov if and only if it is finitely generated and there exists $a,C>0$ such that
$$\alpha_k(\kappa(\gamma))\ge a |\gamma|-C$$
for all $\gamma\in\Gamma$, where $|\gamma|$ is the word length $\gamma$ with respect to some fixed finite
generating set for $\Gamma$.
\end{theorem}

\subsection{Patterson--Sullivan measures} 
For each $k=1,\dots,d-1$, the \emph{$k$-th restricted fundamental weight} is the linear map
\[\omega_k:\mathfrak a\to\Rb\] 
given by $\omega_k(A)=A_1+\dots+A_k$. For any non-empty $\theta\subset\{1,\dots,d-1\}$, let 
\[\mathfrak a_\theta:=\{A\in\mathfrak a:\alpha_k(A)=0\text{ for all }k\notin\theta\}\quad\text{and}\quad\mathfrak a_\theta^+:=\mathfrak a^+\cap\mathfrak a_\theta.\]

Observe that $\{\omega_k|_{\mathfrak a_\theta}:k\in\theta\}$ is a basis for $\mathfrak a_\theta^*$, so we may identify $\mathfrak a_\theta^*$ as the subspace of $\mathfrak a^*$ spanned by $\{\omega_k:k\in\theta\}$. We may then define a projection
\[p_\theta:\mathfrak a\to\mathfrak a_\theta\]
defined by $\omega_k(A)=\omega_k(p_\theta(A))$ for all $k\in\theta$. With this, define the \emph{$\theta$-partial Cartan projection} to be the map
\[\kappa_\theta:=p_\theta\circ\kappa:\PSL(d,\Kb)\to\mathfrak a_\theta^+.\]
Via the identification of $\mathfrak a_\theta^*$ as a subspace of $\mathfrak a^*$, for all $\phi\in\mathfrak a_\theta^*$ and $g\in\PSL(d,\Kb)$, we have
\[\phi(\kappa(g))=\phi(\kappa_\theta(g)).\]

For all $k\in\Delta$, let $\norm{\cdot}_k$ denote the norm on $\wedge^k\Kb^d$ induced by the standard inner product on $\Kb^d$. One may then verify that for all $g\in\PSL(d,\Kb)$ and $x\in\Fc_\Delta$, we have
\[\omega_k(B(g,x))=\log\frac{\norm{g(v_1)\wedge\dots\wedge g(v_k)}_k}{\norm{v_1\wedge\dots\wedge v_k}_k},\]
where $v_1,\dots,v_k$ is some (any) basis of $x^k$. In particular, $\omega_k(B(g,x))$ depends only on $g$ and $x^k$. Thus, for any non-empty $\theta\subset\Delta$, we may define the \emph{$\theta$-partial Iwasawa cocycle}
\[B_\theta:\PSL(d,\Kb)\times\Fc_\theta\to\mathfrak a_\theta\]
by $B_\theta(g,x)=p_\theta(B(g,\sigma(x)))$, where $\sigma:\Fc_\theta\to\Fc_\Delta$ is some (any) right inverse of the forgetful projection $\Fc_\Delta\to\Fc_\theta$. Again, via the identification of $\mathfrak a_\theta^*$ as a subspace of $\mathfrak a^*$, for all $\phi\in\mathfrak a_\theta^*$, $g\in\PSL(d,\Kb)$, and $x\in\Fc_\theta$, we have
\[\phi(B_\theta(g,x))=\phi(B(g,\sigma(x)).\]

Let $\Gamma\subset\PSL(d,\Kb)$ be $\mathsf{P}_\theta$-divergent and let $\phi\in\mathfrak a_\theta^*$. The \emph{$\phi$-critical exponent} of $\Gamma$, denoted $\delta^{\phi}(\Gamma)\in[0,+\infty]$, is the critical exponent of the Poincar\'e series
\[Q^{\phi}_\Gamma(s):=\sum_{\gamma\in\Gamma}e^{-s\phi(\kappa_\theta(\gamma))}=\sum_{\gamma\in\Gamma}e^{-s\phi(\kappa(\gamma))}.\]
If $\delta^\phi(\Gamma)<+\infty$, a \emph{$\phi$-Patterson--Sullivan measure for $\Gamma$} (of dimension $\delta^{\phi}(\Gamma)$) is a probability measure on $\Lambda_\theta(\Gamma)$ such that for all $\gamma\in\Gamma$, $\gamma_*\mu$ and $\mu$ are mutually absolutely continuous and
\[\frac{d\gamma_*\mu}{d\mu}(x)=e^{-\delta^{\phi}(\Gamma)\phi(B_\theta(\gamma^{-1},x))}\]
for $\mu$-almost every $x\in\Lambda_\theta(\Gamma)$. By \cite[Prop.\ 3.2]{CZZ3}, such measures exist. 

In the case when $\Gamma$ is $\mathsf{P}_\theta$-transverse, notice that for any non-empty subset $\theta'\subset\theta$, the forgetful projection $\Fc_\theta\to\Fc_{\theta'}$ restricts to a homeomorphism between $\Lambda_\theta(\Gamma)$ and $\Lambda_{\theta'}(\Gamma)$. Thus, for any $\phi\in\mathfrak a_\theta^*$ with $\delta^\phi(\Gamma)<+\infty$ and any $\phi$-Patterson--Sullivan measure $\mu$ for $\Gamma$ on $\Lambda_\theta(\Gamma)$, we may pushforward $\mu$ by this projection to obtain a measure on $\Lambda_{\theta'}(\Gamma)$. We refer to this pushforward measure also as a \emph{$\phi$-Patterson--Sullivan measure for $\Gamma$ on $\Lambda_{\theta'}(\Gamma)$}. Notice here that $\theta'$ need not be symmetric, and $\phi$ need not lie in $\mathfrak a_{\theta'}^*$.

For every $k\in\{1,\dots,d-1\}$, set 
$$
\theta_k:=\{k-1,k,k+1,d-k-1,d-k,d-k+1 \} \cap \{ 1,\dots, d-1\}.
$$
In the remainder of the paper, when discussing Patterson--Sullivan measures, we will largely specialize to the case when $\Gamma$ is $\mathsf{P}_{\theta_k}$-Anosov and $\mu$ is a $\alpha_k$-Patterson--Sullivan measure for $\Gamma$ on $\Lambda_k(\Gamma)$. (Notice that $\alpha_1=2\omega_1-\omega_2$, $\alpha_{d-1}=2\omega_{d-1}-\omega_{d-2}$ and $\alpha_k=2\omega_k-\omega_{k-1}-\omega_{k+1}$ for every $k\in\{2,\dots,d-2\}$, so we may view $\alpha_k$ as an element in $\mathfrak a_{\theta_k}^*$ via the usual identification of $\mathfrak a_{\theta_k}^*$ as a subspace of $\mathfrak a^*$.) Since $\Gamma$ is $\mathsf{P}_{\theta_k}$-Anosov and $\alpha_k$ is positive on $\mathfrak a_{\theta_k}^+$, it follows from Sambarino \cite[Cor.\ 5.4.3]{sambarino-dichotomy} that $\delta^{\alpha_k}(\Gamma)<+\infty$, and \cite[Thm.\ A]{sambarino-dichotomy} that $\mu$ is unique.

\subsection{Hyperconvex groups}\label{sec: hyperconvex}
Let $k\in\{1,\dots,d-1\}$. With $\theta_k \subset \Delta$ as in the previous subsection, a subgroup $\Gamma\subset\mathsf{PSL}(d,\mathbb K)$ is $(d-k)$-{\em hyperconvex} if it is $\mathsf P_{\theta_k}$-Anosov, and for all mutually distinct triples $x,y,z\in \Lambda_\theta(\Gamma)$ we have
$$\Big((x^k\cap z^{d-k+1})+z^{d-k-1}\Big)\cap \Big( (y^k\cap z^{d-k+1})+z^{d-k-1}\Big)=z^{d-k-1}.$$

Notice that $\Gamma$ is $(d-1)$-hyperconvex if and only if for all mutually distinct triples $x,y,z\in \Lambda_{\theta_1}(\Gamma)$, we have
$$x^1+y^1+z^{d-2}=\mathbb K^d.$$
This condition was called \emph{$(1,1,2)$-hyperconvexity} by 
Pozzetti--Sambarino--Wienhard \cite{PSW1}. They showed that in this setting the $\alpha_1$-Patterson--Sullivan measure behaves much like a conformal measure in
the rank one setting. They also used this observation to bound the $\alpha_1$-critical exponent.

\begin{theorem}[{\cite[Cor.\ 6.9 and 7.3]{PSW1}}]\label{thm: PSW}
If $\Gamma\subset\mathsf{PSL}(d,\mathbb C)$ is $(d-1)$-hyperconvex, then the Hausdorff dimension of
$\Lambda_1(\Gamma)$ is equal to $\delta^{\alpha_1}(\Gamma)$. Moreover, $\delta^{\alpha_1}(\Gamma)\le 2$.
\end{theorem}

Farre--Pozzetti--Viaggi \cite{FPV} showed that hyperconvex groups are always virtually isomorphic to convex cocompact Kleinian groups.
We recall that two groups $G$ and $H$ are \emph{virtually isomorphic} if there exists finite index subgroups $G'\subset G$ and $H'\subset H$ and
finite normal subgroups $M\subset G'$ and $N\subset H'$ so that $G'/M$ is isomorphic to $H'/N$.

\begin{theorem}[{\cite[Thm.\ A]{FPV}}]\label{convex cocompact}
If $\Gamma\subset\mathsf{PSL}(d,\mathbb C)$ is $(d-k)$-hyperconvex for some $k\in\{1,\dots,d-1\}$, then 
$\Gamma$ is virtually isomorphic to a convex cocompact subgroup of $\mathsf{SL}(2,\mathbb C)$.
\end{theorem}

They further adapt Ahlfors' proof \cite{ahlfors}  that a convex cocompact Kleinian group is not cocompact, then its limit set
has measure zero, to control the 2-dimensional Hausdorff measure of limit sets of hyperconvex groups.

\begin{theorem}[{\cite[Thm.\ B]{FPV}}]\label{2-sphere}
If $\Gamma\subset\mathsf{PSL}(d,\mathbb C)$ is $(d-k)$-hyperconvex and $\partial \Gamma$ is not a two-sphere, then the 2-dimensional Hausdorff measure of 
$\Lambda_k(\Gamma)$ is zero. 
\end{theorem}

We will often be able to reduce facts about $(d-k)$-hyperconvex groups to facts about $(d-1)$-hyperconvex groups using the exterior power representation. 
We explain this standard techniques  more carefully in Appendix A, but we state the conclusion of this analysis here.

\begin{proposition}\label{prop: wedge}
If $k\in\{2,\dots,d-1\}$, then there exist $d_k$, a representation $\eta_k:\mathsf{PSL}(d,\Kb)\to\mathsf{PSL}(d_k,\Kb)$ and a smooth embedding
$E_k:\Fc_{k-1,k,k+1}(\Kb^d)\to \mathbb \Fc_{1,2}(\Kb^{d_k})$, which is holomorphic when $\Kb=\Cb$, such that
if  $\Gamma\subset\PSL(d,\Kb)$ is  $(d-k)$-hyperconvex subgroup then the following hold.
\begin{enumerate}
\item The subgroup $\eta_k(\Gamma)\subset\PSL(d_k,\Kb)$ is $(d_k-1)$-hyperconvex, and
$$
\Lambda_{1,2}(\eta_k(\Gamma))=E_k(\Lambda_{k-1,k,k+1}(\Gamma)).
$$
\item $\alpha_1(\kappa(\eta_k(g))) =\alpha_k(\kappa(g))$ for all $g \in \PSL(d,\Kb)$ and hence  $\delta^{\alpha_1}(\eta_k(\Gamma))=\delta^{\alpha_k}(\Gamma)=:\delta$.
\item $E_k$ descends to a smooth embedding $E_k^1 : \Gr_k(\Kb^d) \rightarrow \Pb(\Rb^{d_k})$. 
\item  If $\mu$ is the $\alpha_k$-Patterson-Sullivan measure on $\Lambda_k(\Gamma)$ for $\Gamma$ of dimension $\delta$, 
then $(E_k^1)_*(\mu)$ is the $\alpha_1$-Patterson-Sullivan measure on $\Lambda_1(\eta(\Gamma))$ for $\eta_k(\Gamma)$ of dimension $\delta$.
\end{enumerate}
(In the above, $\alpha_{1}$ is a root of $\PSL(d_k,\Kb)$ while $\alpha_{k}$ is a root of $\PSL(d,\Kb)$.)
\end{proposition}

As an example of how this is applied, notice that Theorem \ref{thm: PSW} and Proposition \ref{prop: wedge} have the following immediate corollary.

\begin{corollary}
If $\Gamma\subset\mathsf{PSL}(d,\mathbb C)$ is $(d-k)$-hyperconvex, then the Hausdorff dimension of
$\Lambda_k(\Gamma)$ is equal to $\delta^{\alpha_k}(\Gamma)$. Moreover, $\delta^{\alpha_k}(\Gamma)\le 2$.
\end{corollary}

\subsection{Kleinian groups} 

Recall that $\PSL(2,\Cb)$ can be identified with $\mathsf{Isom}_0(\Hb^3)$, the connected component of the identity in the isometry group of real hyperbolic 3-space, in such a way that the geodesic boundary $\partial \Hb^3$ identifies with $\mathbb{CP}^1$. Under this identification, 
\begin{enumerate}
\item A subgroup $\Gamma \subset \PSL(2,\Cb)=\mathsf{Isom}_0(\Hb^3)$ is $\Psf_1$-transverse if and only if it is discrete. 
\item A subgroup $\Gamma \subset \PSL(2,\Cb)=\mathsf{Isom}_0(\Hb^3)$ is $\Psf_1$-Anosov if and only if it is convex cocompact as a subgroup of  $\mathsf{Isom}_0(\Hb^3)$. 
\end{enumerate} 

We further note that the unique $\mathsf{PSU}(d)$-invariant probability measure on $\mathbb{CP}^1$ is a $\alpha_1$-Patterson--Sullivan measure of dimension two for any uniform lattice in  $\PSL(2,\Cb)$. 

\begin{remark} There exist finitely generated Kleinian groups with critical exponent 2 whose limit set in $\mathbb{CP}^1$ has Hausdorff
dimension two, but has Hausdorff measure zero (see \cite{canary-laplace}). Recall that all Kleinian groups are $\Psf_1$-transverse. \end{remark}

%% file: ahlfors.tex

\section{Ahlfors regularity}

Given a metric space $X$, a Borel measure $\mu$ on $X$ and $\delta > 0$, we say that $\mu$ is \emph{$\delta$-Ahlfors regular} if there are constants $D\ge 1$ and $t_0>0$ such that for every $x\in X$ and every $0<t<t_0$, we have
\[\frac{1}{D}t^\delta\le \mu (B(x,t))\le Dt^\delta,\]
where $B(x,t)$ is the ball centered at $x$ of radius $t$. Let $k\in\{1,\dots,d-1\}$. In this section we prove that
if $\Gamma$ is $(d-k)$-hyperconvex, then its $\alpha_k$-Patterson--Sulivan measure is Ahlfors regular on
$\Lambda_k(\Gamma)$.

\begin{theorem}\label{prop: Ahlfors regularity}
If $\Gamma\subset \PSL(d,\mathbb C)$ is $(d-k)$-hyperconvex, then the $\alpha_k$-Patterson-Sullivan measure $\mu$ for 
$\Gamma$ on $\Lambda_k(\Gamma)$ is $\delta^{\alpha_k}(\Gamma)$-Ahlfors regular. In particular, $\Lambda_k(\Gamma)$ has finite, non-zero 
$\delta^{\alpha_k}(\Gamma)$-dimensional Hausdorff measure.
\end{theorem}

Theorem \ref{finite Hausdorff measure} follows immediately. Also, as a consequence of Theorem \ref{prop: Ahlfors regularity} and results from Section \ref{sec: hyperconvex}, we may establish the first claim in Theorem \ref{rigidity result}.

\begin{corollary}
\label{first claim}
If $\Gamma\subset \PSL(d,\mathbb C)$ is $(d-k)$-hyperconvex, then $\delta^{\alpha_k}(\Gamma)=2$
if and only if $\Gamma$ is virtually isomorphic to a uniform lattice in $\mathsf{PSL}(2,\mathbb C)$.
\end{corollary}

\begin{proof}
Suppose that $\delta^{\alpha_k}(\Gamma)=2$. By Theorem \ref{prop: Ahlfors regularity}, the $2$-dimensional Hausdorff measure of $\Lambda_k(\Gamma)$ is positive, so Theorem \ref{2-sphere} implies that $\partial\Gamma$ is a $2$-sphere. Then by Theorem \ref{convex cocompact}, $\Gamma$ is virtually isomorphic to a uniform lattice in $\PSL(2,\Cb)$. 

Conversely, if $\Gamma$ is virtually isomorphic to a uniform lattice in $\PSL(2,\Cb)$, then $\partial\Gamma$ is a two-sphere. Therefore, $\Lambda_k(\Gamma)$ has
Hausdorff dimension at least two, so, by Theorem \ref{thm: PSW}, $\delta^{\alpha_1}(\Gamma)\ge 2$. However, Theorem \ref{thm: PSW} also implies
that  $\delta^{\alpha_1}(\Gamma)\le 2$, so  $\delta^{\alpha_1}(\Gamma)= 2$.
\end{proof}

\subsection{Transverse representations and shadows}
Let $\Omega \subset \mathbb P(\Rb^{d_0})$ be a properly convex domain, i.e. is convex and bounded in some affine chart of $\mathbb P(\Rb^{d_0})$. For any $a\in\partial\Omega$, a projective hyperplane $H\subset\Pb(\Rb^{d_0})$ is a \emph{supporting hyperplane to $\Omega$ through $a$} if $a \in H$ and $H\cap\Omega=\emptyset$. We say that $a\in\partial\Omega$ is a \emph{$\Cc^1$ point} if there is a unique supporting hyperplane to $\Omega$ through $a$.

The \emph{Hilbert metric} $\mathrm{d}_\Omega$ on $\Omega$ is defined as follows: Given $p,q\in\Omega$, let $a,b\in\partial\Omega$ be distinct points such that $a,p,q,b$ lie in this order along the projective line segment in $\Omega$ with $a$ and $b$ as endpoints. Then
\[\mathrm{d}_\Omega(p,q):=\log\frac{\abs{q-a}\abs{p-b}}{\abs{p-a}\abs{q-b}},\]
where $\abs\cdot$ denotes some (any) norm on some (any) affine chart of $\Pb(\Rb^{d_0})$ that contains $a$, $p$ ,$q$ and $b$. It is well-known that ${\rm d}_\Omega$ is an $\Aut(\Omega)$-invariant, Finsler metric on $\Omega$, where $\Aut(\Omega)$ denotes the subgroup of $\PGL({d_0},\Rb)$ that leaves $\Omega$ invariant.

Let $\Gamma_0\subset\Aut(\Omega)$ be a discrete subgroup. The \emph{full orbital limit set} of $\Gamma_0$ is 
$$
\Lambda_\Omega(\Gamma_0):=\left\{x\in\partial\Omega\ |\ x=\lim\gamma_n(p)\ \text{ for some }p\in\Omega\text{ and some }\{\gamma_n\}\subset\Gamma_0\right\}.
$$
We may specify a collection of open subsets of $\Lambda_\Omega(\Gamma_0)$ called \emph{shadows} in the following way. For any $p,q\in\Omega$ and $r>0$, let $\Oc_r(p,q)$ denote the set of points $a\in\partial\Omega$ for which the projective line segment $[p,a)$ in $\Omega$ between $p$ and $a$ intersects the open ball $B_\Omega(q,r)$ of radius $r$ (in the Hilbert metric on $\Omega$) centered at $q$. We refer to 
\[\widehat{\Oc}_r(p,q):=\Oc_r(p,q)\cap\Lambda_\Omega(\Gamma_0)\] 
as the \emph{$r$-shadow} in $\Lambda_\Omega(\Gamma_0)$ from $p$ centered at $q$.

Given a properly convex domain $\Omega\subset\Pb(\Rb^{d_0})$, we say that a discrete subgroup $\Gamma_0\subset\Aut(\Omega)$ is \emph{projectively visible} if every point in $\Lambda_\Omega(\Gamma_0)$ is a $\Cc^1$ point, and any two distinct points in $\Lambda_\Omega(\Gamma_0)$ are joined by a projective line segment in $\Omega$. For any such $\Gamma_0$, a representation $\rho : \Gamma_0 \rightarrow \PSL(d,\Kb)$ is called \emph{$\mathsf{P}_\theta$-transverse} if there exists a continuous embedding 
\begin{align*}
\xi : \Lambda_\Omega(\Gamma_0) \rightarrow \mathcal F_\theta(\Kb^d)
\end{align*}
with the following properties: 
\begin{enumerate}
\item $\xi$ is $\rho$-equivariant, i.e. $\rho(\gamma) \circ \xi = \xi \circ \gamma$ for all $\gamma \in \Gamma$, 
\item $\xi$ is transverse, i.e. $\xi(\Lambda_\Omega(\Gamma))$ is a transverse subset of $\mathcal F_\theta(\Kb^d)$, 
\item $\xi$ is strongly dynamics preserving, i.e. if $(\gamma_i)_{i=1}^\infty$ is a sequence in $\Gamma$ so that $\gamma_i(p)\to a\in\Lambda_\Omega(\Gamma)$ and 
$\gamma_i^{-1}(p)\to b\in\Lambda_\Omega(\Gamma)$ for some (any) $p\in\Omega$, then
$\rho(\gamma_i)(x)\to\xi(a)$ if $x\in\mathcal F_\theta(\Kb^d)$  is transverse to $\xi(b)$.
\end{enumerate}
We refer to $\xi$ as the \emph{limit map} of $\rho$ (note that this is unique to $\rho$). 

It follows from Proposition~\ref{prop:characterizing convergence in general} that if $\rho:\Gamma_0\to\PSL(d,\Kb)$ is a $\mathsf P_\theta$-transverse representation with limit map $\xi$, then $\rho(\Gamma_0)$ is a $\mathsf{P}_\theta$-transverse group with $\xi(\Lambda_\Omega(\Gamma_0))$ as its $\theta$-limit set. The following theorem states that the converse is also true when $\theta=\{k,d-k\}$.

\begin{theorem}[{\cite[Thm.\ 6.2]{CZZ3},\cite[Cor.\ 1.30]{zimmer-jdg}}]\label{transverse rep}
If $\Gamma \subset \PSL(d,\Kb)$ is a $\mathsf{P}_{k,d-k}$-transverse subgroup for some $k\le d-k$, then there exist some $d_0\in\Nb$, a properly convex domain $\Omega\subset\Pb(\Rb^{d_0})$, a projectively visible subgroup $\Gamma_0\subset\Aut(\Omega)$, and a faithful 
$\mathsf{P}_{k,d-k}$-transverse representation $\rho: \Gamma_0 \rightarrow \PSL(d, \Kb)$ with limit map 
$\xi:\Lambda_{\Omega}(\Gamma_0)\to\mathcal F_{k,d-k}$ so that 
\begin{enumerate}
\item $\rho(\Gamma_0)=\Gamma$,
\item  $\xi(\Lambda_{\Omega}(\Gamma_0)) = \Lambda_{k,d-k}(\Gamma)$,
\item $\alpha_1(\kappa(\gamma))=\alpha_k(\kappa(\rho(\gamma)))$ for all $\gamma\in\Gamma_0$.
\end{enumerate} 
Furthermore, if $\Gamma$ is $\Psf_{k,d-k}$-Anosov, then $\Gamma_0$ acts cocompactly on the convex hull in $\Omega$ of $\Lambda_\Omega(\Gamma_0)$.
\end{theorem}

If $\Gamma\subset\PSL(d,\Kb)$ is a $\mathsf{P}_{1,d-1}$-transverse subgroup, then by making a choice of $d_0$, $\Omega$, $\Gamma_0$, and $\rho$ as given by Theorem \ref{transverse rep} (with $k=1$), we may define shadows in $\Lambda_1(\Gamma)$ as the images of the shadows in $\Lambda_\Omega(\Gamma_0)$ under $\xi^1$, where $\xi=(\xi^1,\xi^{d-1})$ is the limit map of $\rho$. In our previous work \cite[Prop.\ 7.1]{CZZ3}, we proved a shadow lemma for these shadows, i.e. we estimated the $\phi$-Patterson--Sullivan measures of the shadows in $\Lambda_1(\Gamma)$. Also, assuming an analog of the $(d-1)$-hyperconvexity assumption for transverse representations, we estimated the in-radii and out-radii of these shadows \cite[Thm.\ 5.1 and 7.1]{CZZ2}. We will now state a specialization of both of these results to the setting where $\Gamma$ is $(d-1)$-hyperconvex and $\phi=\alpha_1$.

\begin{theorem}[{\cite[Prop.\ 7.1]{CZZ3}, \cite[Thm.\ 5.1 and 7.1]{CZZ2}}]\label{prop:shadow estimates} Let $\Gamma\subset\PSL(d,\Kb)$ be a $(d-1)$-hyperconvex subgroup and let $\mu$ be the $\alpha_1$-Patterson--Sullivan measure for $\Gamma$ on $\Lambda_1(\Gamma)$. Let $d_0\in\Nb$, $\Omega\subset\Pb(\Rb^{d_0})$ be a properly convex domain, $\Gamma_0\subset \Aut(\Omega)$ be a projectively visible group, and $\rho:\Gamma_0\to\PSL(d,\Kb)$ be a faithful, $\mathsf{P}_{1,d-1}$-transverse representation with limit map $\xi:\Lambda_\Omega(\Gamma_0)\to\Fc_{1,d-1}$ satisfying properties (1)--(3) of Theorem \ref{transverse rep} (with $k=1$). Fix a base point $o\in\Omega$. There exists $r_0>0$ such that for all $r>r_0$, there exists $C(r)>1$ for which the following hold:
\begin{enumerate}
\item[(i)] For all $\gamma \in \Gamma_0$,
\begin{align*}
\frac{1}{C(r)}e^{-\delta \alpha_1(\kappa(\rho(\gamma)))} \le \mu\Big(\xi^1\big(\widehat{\Oc}_r(o,\gamma(o)) \big) \Big)  \le C(r)e^{-\delta\alpha_1(\kappa(\rho(\gamma)))}.
\end{align*}
\item[(ii)] If $a\in\Lambda_\Omega(\Gamma_0)$, $p$ lies in the projective line segment $[o,a)$ in $\Omega$ with endpoints $o$ and $a$, and $\gamma\in\Gamma_0$ satisfies ${\rm d}_\Omega(p,\gamma(o))<r$, then
\begin{align*}
B\left(\xi^1(a), \frac{1}{C(r)} e^{-\alpha_1(\kappa(\rho(\gamma)))}\right)\cap \Lambda_1(\Gamma)\subset\xi^1\left(\widehat{\mathcal O}_r(o,p)\right)   \subset B\left(\xi^1(a),C(r)e^{-\alpha_1(\kappa(\rho(\gamma)))}\right)\cap \Lambda_1(\Gamma),
\end{align*}
where $B(x,t)$ denotes the ball centered at $x\in\Fc_1(\Kb^d)$ of radius $t$ in the Riemannian metric on $\Fc_1(\Kb^d)$.
\end{enumerate}
\end{theorem} 

The results in~\cite{CZZ2} used closed shadows, but due to the nesting properties of these sets, the difference is immaterial.

\subsection{Proof of Theorem \ref{prop: Ahlfors regularity} } Using Proposition~\ref{prop: wedge} it suffices to consider the case when $k=1$.

Suppose that $\Gamma\subset \PSL(d,\mathbb K)$ is $(d-1)$-hyperconvex and $\mu$ is the  $\alpha_1$-Patterson--Sullivan measure  
for $\Gamma$.
By Theorem \ref{transverse rep}, there exist  $d_0\in\Nb$, a properly convex domain $\Omega\subset\Pb(\Rb^{d_0})$, 
a projectively visible subgroup $\Gamma_0\subset\Aut(\Omega)$, and a faithful, $\mathsf{P}_{1,d-1}$-transverse representation 
$\rho:\Gamma_0\to\PSL(d,\Kb)$ with image $\Gamma$ and limit map $\xi:\Lambda_\Omega(\Gamma_0)\to\Lambda_{1,d-1}(\Gamma)$ 
that satisfy properties (1)--(3) of Theorem \ref{transverse rep} (with $k=1$).
Since $\Gamma$ is $\mathsf{P}_{1,d-1}$-Anosov, $\Gamma_0$ acts cocompactly on the convex
hull $\mathcal C$ in $\Omega$ of $\Lambda_\Omega(\Gamma_0)$.  So, if we choose $o\in \mathcal C$, 
then there exists $R>0$ so that if $p\in\mathcal C$ then there exists $\gamma\in\Gamma_0$ so that
$d_\Omega(p,\gamma(o))<R$ (in the Hilbert metric on $\Omega$). 
We may also assume that $R>r_0$, where $r_0>0$ is the constant in Theorem \ref{prop:shadow estimates}. 

Fix $x\in\Lambda_1(\Gamma)$. Then there exists $a\in\Lambda_\Omega(\Gamma_0)$ such that $\xi^1(a)=x$. 
Fix a base point $o\in\Omega$, and let $[o,a)$ denote the projective line segment in $\Omega$ with endpoints $o$ and $a$. 
Recall that ${\rm d}_\Omega$ denotes the Hilbert metric on $\Omega$ and for any $p\in\Omega$ and $r>0$, 
$B_\Omega(p,r)$ denotes the open Hilbert ball of radius $r$ centered at $p$. 

Iteratively construct a sequence $\{\gamma_n\}$ in $\Gamma_0$ as follows. Set $\gamma_1:=\id$. Then given $\gamma_n$, choose $\gamma_{n+1}\in\{\gamma\in\Gamma_0:B_\Omega(\gamma(o),R)\cap[o,a)\neq\emptyset\}$ such that the endpoint of the interval $B_\Omega(\gamma_n(o),R)\cap[o,a)$ that is furthest from $o$ lies in $B_\Omega(\gamma_{n+1}(o),R)$. Observe that 
${\rm d}_\Omega(\gamma_n(o),\gamma_{n+1}(o))< 2R$ for all integers $n\ge 1$, and 
$[o,a)\subset\bigcup_{n=1}^\infty B_\Omega(\gamma_n(o),R)$.

For each positive integer $n$, set $\beta_n:=\gamma_n^{-1}\gamma_{n+1}$. It is a standard linear algebra fact that for all $g_1,g_2\in\PSL(d,\Kb)$ and all $i\in\{1,\dots,d\}$, we have 
\[\sigma_i(g_1)\sigma_d(g_2)\le\sigma_i(g_1g_2)\le\sigma_i(g_1)\sigma_1(g_2).\]
Thus, for all positive integers $n$, we have 
\[\alpha_1(\kappa(\rho(\gamma_n)))-\log\frac{\sigma_1(\rho(\beta_n))}{\sigma_d(\rho(\beta_n))}\le\alpha_1(\kappa(\rho(\gamma_{n+1})))\le\alpha_1(\kappa(\rho(\gamma_n)))+\log\frac{\sigma_1(\rho(\beta_n))}{\sigma_d(\rho(\beta_n))}.\]
Since ${\rm d}_\Omega(b,\beta_n(b))\le 2R$ for all $n$, if we set 
$$S:=\max\left\{\frac{\sigma_1(\rho(\beta))}{\sigma_d(\rho(\beta))} : \beta\in\Gamma_0\text{ and } {\rm d}_\Omega(b,\beta(b))\le 2R\right\},$$
then for all $n$, we have
\begin{equation}\label{eqn: sequence}\abs{\alpha_1(\kappa(\rho(\gamma_n)))-\alpha_1(\kappa(\rho(\gamma_{n+1})))}\le\log S.\end{equation}

For each $r>r_0$, let $C(r)\ge 1$ be the constant given by Theorem \ref{prop:shadow estimates}. Set $\delta:=\delta^{\alpha_1}(\Gamma)$, $t_0:=\frac{1}{SC(R)}$, and $D:=\max\{C(2R) C(R)^\delta S^\delta,C(R)C(2R)^\delta S^{\delta}\}$. To prove Ahlfors regularity, it suffices to show that for all $0<t<t_0$, we have
\begin{equation}
\label{ahlfors reg}
\frac{1}{D}t^\delta\le \mu (B(x,t))\le Dt^\delta.
\end{equation}

First, we prove the required upper bound for $\mu (B(x,t))$. Since $\rho$ is $\mathsf{P}_{1,d-1}$-divergent and $\{\gamma_n\}$ is a diverging sequence in $\Gamma_0$, we have $\alpha_1(\kappa(\rho(\gamma_n)))\to\infty$, so we may let $\ell$ be the smallest positive integer such that $\alpha_1(\kappa(\rho(\gamma_\ell)))\ge -\log(SC(R)t)$. Since $\gamma_1=\id$ and $0<-\log(SC(R)t)$, observe that $\ell\ge 2$ and $\alpha_1(\kappa(\rho(\gamma_{\ell-1})))<-\log(SC(R)t)$. It now follows from \eqref{eqn: sequence} that
\[-\log(SC(R)t)\le \alpha_1(\kappa(\rho(\gamma_{\ell})))\le -\log(SC(R)t)+\log S=-\log(C(R)t),\]
or equivalently, that
\begin{equation}\label{eqn: Ahlfors 1}\frac{1}{SC(R)}e^{-\alpha_1(\kappa(\rho(\gamma_{\ell})))}\le t\le \frac{1}{C(R)}e^{-\alpha_1(\kappa(\rho(\gamma_{\ell})))}.\end{equation}
Choose $p\in[o,a)$ so that ${\rm d}_\Omega(p,\gamma_{\ell}(o))<R$. By the first inclusion of Theorem \ref{prop:shadow estimates} part (ii) and the second inequality in \eqref{eqn: Ahlfors 1}, we have that
$$B(x,t)\cap\Lambda_1(\Gamma)\subset \xi^1(\widehat{\mathcal O}_R(o,a))\subset \xi^1(\widehat{\mathcal O}_{2R}(o,\gamma_{\ell}(o))).$$
Then the second inequality of Theorem \ref{prop:shadow estimates} part (i) and the first inequality in \eqref{eqn: Ahlfors 1} imply that
$$\mu(B(x,t))\le \mu\left( \xi^1(\widehat{\mathcal O}_{2R}(o,\gamma_{\ell}(o)))\right)\le C(2R)e^{-\delta\alpha_1(\kappa(\rho(\gamma_\ell)))}\le C(2R) (SC(R))^\delta t^{\delta}\le Dt^\delta.$$

The proof of the required lower bound for $\mu (B(x,t))$ is similar. Observe that since $0<t< SC(2R)$ (because $t_0<1\le SC(2R)$) and $\gamma_1=\id$, we may let $m$ be the largest positive integer such that $\alpha_1(\kappa(\rho(\gamma_{m})))\le -\log\frac{t}{SC(2R)}$. Then $\alpha_1(\kappa(\rho(\gamma_{m+1})))\ge-\log\frac{t}{SC(2R)}$, so \eqref{eqn: sequence} implies that
\[-\log \frac{t}{C(2R)}=-\log\frac{t}{SC(2R)}-\log S\le \alpha_1(\kappa(\rho(\gamma_{m})))\le -\log\frac{t}{SC(2R)},\]
or equivalently that
\begin{equation}\label{eqn: Ahlfors2}C(2R)e^{-\alpha_1(\kappa(\rho(\gamma_{m})))}\le t\le SC(2R)e^{-\alpha_1(\kappa(\rho(\gamma_{m})))}.\end{equation}
Choose $p\in[o,a)$ so that ${\rm d}_\Omega(p,\gamma_{m}(o))<R$.  By the second inclusion of Theorem \ref{prop:shadow estimates} part (ii) and the first inequality of \eqref{eqn: Ahlfors2},
$$\xi^1\Big(\widehat{\mathcal O}_R(o,\gamma_{m}(o))\Big)\subset \xi^1\Big(\widehat{\mathcal O}_{2R}(o,p)\Big)\subset B(x, t)\cap\Lambda_1(\Gamma).$$
Then the first inequality of Theorem \ref{prop:shadow estimates} part (i) and the second inequality in \eqref{eqn: Ahlfors2} imply that
\[\mu\Big(B(x,t)\Big)\ge \mu\left( \xi^1(\widehat{\mathcal O}_{R}(o,\gamma_{m}(o)))\right)\ge \frac{1}{C(R)}e^{-\delta\alpha_1(\kappa(\rho(\gamma_{m})))}\ge  \frac{t^\delta}
{C(R)(C(2R)S)^{\delta}}\ge\frac{1}{D}t^\delta.\]
This completes the proof of Equation (\ref{ahlfors reg}) and hence the proof that $\mu$ is $\delta$-Ahlfors regular.

We now recall the standard argument that the support of a $\delta$-Ahlfors regular measure has finite non-zero
$\delta$-dimensional Hausdorff measure.
First, suppose that $\{B(x_i,t_i)\}$ is a countable cover of $\Lambda_1(\Gamma)$ by balls of
radius $t_i<t_0$ based at points $x_i\in\Lambda_1(\Gamma).$
Then, by Equation (\ref{ahlfors reg}) and the fact that $\mu$ is a probability measure, we have
\[\sum_{i=1}^\infty t_i^\delta\ge \frac{1}{D}\sum_{i=1}^\infty \mu(B(x_i,t_i))\ge \frac{1}{D}.\]
It follows that the $\delta$-dimensional Hausdorff measure of $\Lambda_1(\Gamma)$ is
at least $\frac{1}{D}$, hence non-zero.

On the other hand, if $0<t<t_0$, the Vitali covering lemma implies that we may find a finite cover  $\{B(x_i,t)\}_{i=1}^k$ of $\Lambda_1(\Gamma)$
by balls of radius $t$ so that the collection $\{B(x_i,\frac{t}{3})\}_{i=1}^k$ is a disjoint collection of $\frac{t}{3}$-balls.
So
$$1\ge \sum_{i=1}^k \mu\Big(B(x_i,\frac{t}{3})\Big)\ge \frac{k t^\delta}{D3^\delta}$$
which implies that   $\delta$-dimensional Hausdorff measure of $\Lambda_1(\Gamma)$ is
at most $D3^\delta$, hence finite.
\qed

%% file: lattice.tex

\section{Rigidity} 

In this section, we complete the proof of Theorem \ref{rigidity result} by establishing our rigidity result in the lattice case. Recall that a discrete subgroup $\Gamma\subset \PSL(d,\Kb)$ is \emph{strongly irreducible} if every finite index subgroup of $\Gamma$ is irreducible. As before, let $k\in\{1,\dots,d-1\}$.

\begin{theorem}\label{thm:entropy rigidity strongly irred}
If $\Gamma\subset \PSL(d,\mathbb C)$ is a strongly irreducible $(d-k)$-hyperconvex subgroup and 
$$
\delta^{\alpha_k}(\Gamma) = 2,
$$
then $\Gamma$ is a uniform lattice in the image of an irreducible representation $\PSL(2,\Cb) \rightarrow \PSL(d,\Cb)$. 
\end{theorem}

In the case when the Zariski closure of $\Gamma$ is simple, we originally had a proof of Theorem~\ref{thm:entropy rigidity strongly irred} that was more dynamical in flavor (a sketch of this argument is provided at the end of the section). We thank Andres Sambarino for suggesting that we
use Chow's theorem \cite{chow} instead.

Recall that a submanifold in a complex manifold is a \emph{complex submanifold} if it is locally the graph of a holomorphic function. The following lemma is a main ingredient that we will need to apply Chow's theorem.

\begin{lemma}
\label{differentiable general}
If $\Gamma\subset\PSL(d,\mathbb C)$ is $(d-k)$-hyperconvex subgroup and $\partial\Gamma$ is a two-sphere, then $\Lambda_{k-1,k,k+1}(\Gamma)$
is a complex submanifold of $\Fc_{k-1,k,k+1}(\Cb^d)$ (where we use the convention $\Fc_{0,1,2}(\Cb^d) = \Fc_{1,2}(\Cb^d)$ and $\Fc_{d-2,d-1,d}(\Cb^d)=\Fc_{d-2,d-1}(\Cb^d)$). 
\end{lemma}

Assuming Lemma \ref{differentiable general}, we will now prove Theorem \ref{thm:entropy rigidity strongly irred}. Recall that by Schur's lemma, any irreducible Lie subgroup of $\SL(d,\Cb)$ is semisimple, and its center lies in the center of $\SL(d,\Cb)$. Therefore, any irreducible Lie subgroup of $\PSL(d,\Cb)$ is semisimple with trivial center.

We use the following standard fact about irreducible representations of products of semisimple Lie groups. 

\begin{lemma}\label{lem:tensor prod in disguise} Suppose $\mathsf{H}_1, \mathsf{H}_2 \subset \SL(d,\Kb)$ are commuting connected semisimple Lie groups and $\mathsf{H}:=\mathsf{H}_1\mathsf{H}_2$ acts irreducibly on $\Kb^d$. Then there exist irreducible representations $\rho_j : \mathsf{H}_j \rightarrow \SL(d_j, \Kb)$ and a linear identification 
$$
\Kb^d = \Kb^{d_1} \otimes \Kb^{d_2}
$$
such that elements $h = h_1h_2 \in \mathsf{H}$ act by 
$$
h (v_1 \otimes v_2) = (\rho_1(h_1)(v_1)) \otimes (\rho_2(h_2)(v_2)).
$$
\end{lemma} 

\begin{proof} Since $\mathsf{H}_1$ is semisimple, we may decompose its action on $\Kb^d$ into its irreducible factors, i.e. there exist $m_1,\dots, m_k, n_1,\dots, n_k \in \Nb$ and irreducible representations $\tau_j : \mathsf{H}_1 \rightarrow \SL(m_j, \Kb)$ so that we can identify
$$
\Kb^d = \oplus_{j=1}^k (\Kb^{m_j})^{n_j}
$$
and $\mathsf{H}_1$ acts by $(\tau_1^{n_1}, \dots, \tau_k^{n_k})$. 
(See, for example, Serre \cite[Chap. VI.3]{serre}.)
Since $\mathsf{H}_1$ and $\mathsf{H}_2$ commute, $\mathsf{H}_2$ preserves each $(\Kb^{m_j})^{n_j}\subset\Kb^d$. Then, since $\mathsf{H}$ acts irreducibly, we must have $k=1$. 

Relative to the decomposition $\Kb^d = (\Kb^{m_1})^{n_1}$, we can write each element $h \in \mathsf{H}_2$ as 
$$
h = (T_{h,u,v})_{1 \le u, v \le n_1} 
$$
where $T_{h,u,v} : \Kb^{m_1} \rightarrow \Kb^{m_1}$ is linear. Since $\mathsf{H}_1$ and $\mathsf{H}_2$ commute, we have 
$$
T_{h,u,v}  \tau_1(g) =\tau_1(g) T_{h,u,v} 
$$
for all $g \in  \mathsf{H}_1$ and $1 \le u, v \le n_1$. Thus by Schur's lemma, $T_{h,u,v}  = \lambda_{h,u,v}  {\rm Id}_{\Kb^{m_1}}$ for some $\lambda_{h,u,v} \in \Kb$. Define a representation $\rho_2 : \mathsf{H}_2 \rightarrow \SL(n_1, \Kb)$ by $\rho_2(h) = ( \lambda_{h,u,v})_{1 \le u, v \le n_1} $. 

Next we identify 
$$
(\Kb^{m_1})^{n_1} = \Kb^{m_1} \otimes \Kb^{n_1}
$$
by $(v_1,\dots, v_{n_1}) \mapsto v_1 \otimes e_1 + \cdots + v_{n_1} \otimes e_{n_1}$. Under this identification, $\mathsf{H}_1$ acts by 
$$
h_1 (v \otimes w) = (\tau_1(h_1)(v)) \otimes w
$$
and 
$\mathsf{H}_2$ acts by 
$$
h_2 (v \otimes w) = v \otimes (\rho_2(h_2)(w)). 
$$

Notice that if $V \subset \Kb^{n_1}$ is a $\rho_2$-invariant linear subspace, then $\Kb^{m_1} \otimes V\subset \Kb^{m_1}\otimes\Kb^{m_2}$ is $\mathsf{H}$-invariant. Hence $\rho_2$ is irreducible. Finally, setting $d_1 = m_1$, $d_2 = n_1$, and $\rho_1 = \tau_1$ establishes the lemma. 
\end{proof}

\begin{proof}[Proof of Theorem \ref{thm:entropy rigidity strongly irred}]
Let $\mathsf{G}\subset\PSL(d,\Cb)$ denote the real Zariski closure of $\Gamma$ and let $\Gsf^0$ denote the connected component of the identity in $\Gsf$. Since $\Gamma$ is strongly irreducible, it follows that $\Gamma':=\Gamma\cap\Gsf^0$, and hence also $\Gsf^0$, is irreducible. 

Since $\delta^{\alpha_k}(\Gamma) = 2$, by Corollary \ref{first claim}, $\Gamma$ is virtually isomorphic to a uniform lattice in $\mathsf{PSL}(2,\mathbb C)$, so $\partial\Gamma'=\partial\Gamma$ is a two-sphere. We may thus apply Lemma \ref{differentiable general} to deduce that $\Lambda_{k-1,k,k+1}(\Gamma)$ is a complex submanifold in $\Fc_{k-1,k,k+1}(\Cb^d)$, which is necessarily biholomorphic to $\mathbb{CP}^1$. 
We thus have an isomorphism of biholomorphism groups 
$$
\mathsf{Aut}_{hol}(\Lambda_{k-1,k,k+1}(\Gamma)) \cong \mathsf{Aut}_{hol}(\mathbb{CP}^1)=\PSL(2,\Cb).
$$
In what follows we identify 
$$
\mathsf{Aut}_{hol}(\Lambda_{k-1,k,k+1}(\Gamma)) =\PSL(2,\Cb).
$$

By Chow's theorem \cite{chow}, $\Lambda_{k-1,k,k+1}(\Gamma)$ is a complex (and also real) algebraic subvariety of $\Fc_{k-1,k,k+1}(\Cb^d)$. Since $\Gsf$ is the real Zariski closure of $\Gamma$ and $\Lambda_{k-1,k,k+1}(\Gamma)$ is $\Gamma$-invariant, it follows that $\Gsf$ preserves $\Lambda_{k-1,k,k+1}(\Gamma)$. Since $\mathsf{G} \subset \PSL(d,\Cb)$ acts on $\Fc_{k-1,k,k+1}(\Cb^d)$, and hence on $\Lambda_{k-1,k,k+1}(\Gamma)$, by biholomorphisms, this $\Gsf$-action defines a homomorphism 
$$
q_k : \Gsf \rightarrow \mathsf{Aut}_{hol}(\Lambda_{k-1,k,k+1}(\Gamma)) = \PSL(2,\Cb).
$$

\medskip

\textbf{Claim:} $q_k(\Gsf^0) = \PSL(2,\Cb)$. 

\medskip

\emph{Proof of Claim:} Since $\Gamma'$ acts minimally and as a convergence group on $\Lambda_{k-1,k,k+1}(\Gamma)$, the same holds for the $q_k(\Gamma')$-action on $\mathbb{CP}^1$. The minimality of this action implies that $q_k(\Gsf^0)$ is irreducible, and hence semisimple. The fact that this is a convergence group action implies that $q_k(\Gsf^0)$ is non-compact. Since $q_k(\Gsf^0)$  is connected, semisimple, and non-compact, it follows that $q_k(\Gsf^0)$ is either equal to $\PSL(2,\Cb)$ or conjugate to $\PSL(2,\Rb)$. Since $q_k(\Gamma')$ acts minimally on $\mathbb{CP}^1$, we must have $q_k(\Gsf^0)=\PSL(2,\Cb)$. \hfill $\blacktriangleleft$

\medskip

Since $\Gsf^0$ is irreducible,  we observed earlier that it is semisimple and its center is trivial, so we can decompose it into a product of simple Lie groups.  Fix a simple factor $\Gsf_1$ such that $q_k(\Gsf_1)$ is non-trivial. Then $q_k(\Gsf_1)$ is normal in $q_k(\Gsf^0)=\PSL(2,\Cb)$ and $\PSL(2,\Cb)$ is simple, so $q_k(\Gsf_1)=\PSL(2,\Cb)$. Further, since $\Gsf_1$ is simple and has trivial center (since $\Gsf^0$ has trivial center), we see that $q_k|_{\Gsf_1}:\Gsf_1\to\PSL(2,\Cb)$ is an isomorphism. Since $q_k(\Gsf_1)=\PSL(2,\Cb)$ commutes with the images of the other simple factors and $\PSL(2,\Cb)$ has trivial center, $\ker q_k$ coincides with the product of the other simple factors. Hence we can write $\Gsf^0$ as the direct product
\[\Gsf^0=\Gsf_1\times\Gsf_2\]
where $\Gsf_2=\ker q_k$.
We will show that $\Gsf_2=1$, which will imply that $q_k : \Gsf^0 \rightarrow \PSL(2,\Cb)$ is an isomorphism. 

For $j=1,2$, let $\Gsf_j'$ denote the preimage of $\Gsf_j$ under the quotient map $\SL(d,\Cb)\rightarrow \PSL(d,\Cb)$. Notice that $\Gsf_1'$ and $\Gsf_2'$ are commuting subgroups of $\SL(d,\Cb)$, so by Lemma~\ref{lem:tensor prod in disguise}, there exist integers $d_1$ and $d_2$, representations $\rho_1 :\Gsf_1' \rightarrow \SL(d_1,\Cb)$ and $\rho_2: \Gsf_2' \rightarrow \SL(d_2, \Cb)$ and an identification 
$$
\Cb^d = \Cb^{d_1} \otimes \Cb^{d_2}
$$
such that an element $g =g_1g_2 \in \Gsf_1'\Gsf_2'$ acts by 
$$
g(v \otimes w)= (\rho_1(g_1)(v)) \otimes (\rho_2(g_2)(w)).
$$

\medskip

\textbf{Claim:} If $\ell \in \{k-1,k,k+1\}\cap\{1,\dots,d-1\}$, then $$
\Lambda_\ell(\Gamma) \subset \{ V \otimes \Cb^{d_2}  : V \in \Gr_{\ell/d_2}(\Cb^{d_1}) \}. 
$$

\medskip

\emph{Proof of Claim:} Let $\gamma \in \Gamma'$ be an infinite order element and let $U_\gamma^+ \in \Lambda_\ell(\Gamma)$ be the attracting fixed point of $\gamma$. Then $U_\gamma^+$ is a direct sum of generalized eigenspaces of $\gamma$ (corresponding to the largest $\ell$ eigenvalues). If we write $\gamma =\gamma_1\gamma_2$ where $\gamma_j \in \Gsf_j$, then 
$$
U_\gamma^+ = \oplus_i V_i \otimes W_i
$$
where each $V_i$ is a generalized eigenspace of $\rho_1(\gamma_1)$ and each $W_i$ is a generalized eigenspace of $\rho_2(\gamma_2)$.  Since  $\Gsf_2 = \ker q_k$,  for all $g_2 \in \Gsf_2$, we have $g_2(U_\gamma^+)= U_\gamma^+$, so  the union of subspaces $\cup_iW_i$ is invariant under $\rho_2(g_2)$. Since $\rho_2$ is strongly irreducible $W_i=\Cb^{d_2}$ for all $i$. Then 
$$
U_\gamma^+ = V \otimes \Cb^{d_2}
$$
where $V = \sum_i V_i$. Since the set of attracting fixed points of infinite order elements in $\Gamma$ is dense in $\Lambda_\ell(\Gamma)$ and $\Gamma$ acts minimally on $\Lambda_\ell(\Gamma)$, the claim follows. \hfill $\blacktriangleleft$

\medskip

The above claim implies that there is some $\ell\in\{k-1,k\}$ such that $\frac{\ell}{d_2}$ and $\frac{\ell+1}{d_2}$ are both integers, so we must have $d_2 = 1$. This in turn implies that $\Gsf_2$ is a finite, normal subgroup of $\Gsf^0$. Since $\Gsf^0$ is connected, its action by conjugation on $\Gsf_2$ is necessarily trivial, so $\Gsf_2$ must be in the center of $\Gsf^0$, which is trivial. Thus, $\Gsf^0=\Gsf_1$ and $q_k : \Gsf^0 \rightarrow \PSL(2,\Cb)$ is an isomorphism.

To finish the proof, we need to show that $\Gsf = \Gsf^0$. Since $\Gsf^0\subset\PSL(d,\Cb)$ is irreducible and isomorphic to $\PSL(2,\Cb)$, we may view it as the image of an injective irreducible representation $\tau:\PSL(2,\Cb)\to\PSL(d,\Cb)$. By the classification of irreducible representations of $\PSL(2,\Cb)$, there is a linear identification between $\Cb^d$ and the symmetric tensor ${\rm Sym}^{d-1}(\Cb^2)$ with the defining property that
\[\tau(g)[v_1\odot\dots\odot v_{d-1}]=[g(v_1)\odot\dots\odot g(v_{d-1})]\]
for all $v_1,\dots,v_{d-1}\in\Cb^2$ and all $g\in\PSL(2,\Cb)$. The standard (Hermitian) inner product $\langle\cdot,\cdot\rangle$ on $\Cb^2$ induces an inner product $(\cdot,\cdot)$ on ${\rm Sym}^{d-1}(\Cb^2)$ with the defining property 
\[(v_1\odot\dots\odot v_{d-1},u_1\odot\dots\odot u_{d-1}):=\sum_{\sigma\in S_{d-1}}\prod_{i=1}^{d-1}\langle v_i,u_{\sigma(i)}\rangle,\]
where $S_{d-1}$ is the symmetric group on the set $\{1,\dots,d-1\}$. It is straightforward to verify that for all $g\in\PSL(2,\Cb)$ and $j\in\{1,\dots,d-1\}$, 
\[\alpha_j(\kappa(\tau(g)))=\alpha_1(\kappa(g)),\]
where $\kappa$ on the left and the right are the Cartan projections for $\PSL(d,\Cb)$ and $\PSL(2,\Cb)$ with respect to $(\cdot,\cdot)$ and $\langle\cdot,\cdot\rangle$ respectively. In particular, 
\[\alpha_j(\kappa(\gamma))=\alpha_k(\kappa(\gamma))\]
for all $\gamma \in \Gamma'$ and $j,k\in\{1,\dots,d-1\}$. It now follows from Theorem \ref{Anosov singular values} that $\Gamma'$, and hence $\Gamma$, is a $\Psf_\Delta$-Anosov subgroup of $\PSL(d,\Cb)$.

The map 
$$
[v] \in \mathbb{CP}^1 \mapsto [v \odot \cdots \odot v] \in \Pb( {\rm Sym}^{d-1}(\Cb^2))
$$
is a $\tau$-equivariant embedding that is a biholomorphism onto its image, which is the limit set $\Lambda_1(\Gamma)=\Lambda_1(\Gamma')$. Further, $\Gsf$ acts by biholomorphisms on $\Lambda_1(\Gamma)$, so we have a homomorphism
$$
q_1 : \Gsf \rightarrow \mathsf{Aut}_{hol}(\Lambda_1(\Gamma))
$$
from $\Gsf$ to the biholomorphism group of $\Lambda_1(\Gamma)$. Since $\Lambda_1(\Gamma)$ contains $(d+1)$-points in general position (it is the image of the Veronese embedding from $\mathbb{CP}^1$ to $\mathbb{CP}^{d-1}$), the only element of $\Gsf$ acting trivially on $\Lambda_1(\Gamma)$ is the identity. 
Hence $\ker q_1$ is trivial. At the same time, since $\Gsf^0$ is isomorphic to $\PSL(2,\Cb)\cong \mathsf{Aut}_{hol}(\Lambda_1(\Gamma))$, $q_1|_{\Gsf^0}$ is an isomorphism. Thus $\Gsf =\Gsf^0$. 
\end{proof}

It now remains to prove Lemma \ref{differentiable general}. To do so, we first observe from the following result of Pozzetti--Sambarino--Wienhard \cite{PSW1}  (see also Zhang--Zimmer \cite{ZZ} for the case when $\Kb=\Rb$) that $\Lambda_{1,2}(\Gamma)$ is a differentiable submanifold when 
$\Gamma$ is $(d-1)$-hyperconvex.

\begin{theorem}{\rm(\cite[Prop.\ 7.1]{PSW1})}\label{local convergence}
Suppose that $\Gamma\subset\PSL(d,\Kb)$ is $(d-1)$-hyperconvex subgroup. If $\{(z_n,w_n)\}$ is a sequence of distinct points in $\Lambda_{\theta_1}(\Gamma)$ that converges to $(x,x)$, then
$$z_n^1\oplus y_n^1\to x^2.$$
\end{theorem}

Thus, in the case that $\Kb=\mathbb C$ and $\partial\Gamma$ is a two-sphere, Theorem~\ref{local convergence} immediately implies the following consequence. (For comparison, see
\cite[Prop.\ 7.4]{PSW1} and \cite[Thm.\ 1.1]{ZZ}.)

\begin{corollary}
\label{differentiable k=1}
If $\Gamma\subset\PSL(d,\mathbb C)$ is $(d-1)$-hyperconvex subgroup and $\partial\Gamma$ is a two-sphere, then $\Lambda_1(\Gamma)$
is a $\mathcal C^1$-submanifold of $\mathbb P(\mathbb C^d)$. Furthermore, for every $x\in\Lambda_{1,2}(\Gamma)$, we have
 \begin{equation}\label{eqn:the tangent spaces}
 T_{x^1}\Lambda_1(\Gamma) = T_{x^1}\Pb(x^2).
 \end{equation}
\end{corollary}

In the above corollary, $T_{x^1}\Pb(x^2)$ is a complex linear subspace of $T_{x^1} \Pb(\Cb^d)$. So to further upgrade $\Lambda_1(\Gamma)$ to a complex submanifold of $\Pb(\Cb^d)$, we use the following general statement.

\begin{lemma} \label{complex submanifold}
Suppose $N$ is a $\Cc^1$-submanifold of a complex manifold $M$. If for every $x \in N$ the tangent space $T_x N$ is a complex linear subspace of $T_x M$, then $N$ is a complex submanifold. In particular, if $\Gamma\subset\PSL(d,\mathbb C)$ is $(d-1)$-hyperconvex subgroup and $\partial\Gamma$ is a two-sphere, then $\Lambda_1(\Gamma)$
is a complex submanifold of $\mathbb P(\mathbb C^d)$. 
\end{lemma} 

\begin{proof} 
Since being a complex manifold is a local condition, it suffices to prove the following: If $U \subset \Cb^n$ is an open connected set and $f : U \rightarrow \Cb^m$ is a $\Cc^1$-smooth function such that every tangent plane to the $\Cc^1$-submanifold
$$
N:=\{ (z,f(z))\in\Cb^{m+n} : z \in U\}
$$
is complex linear, then $f$ is holomorphic. Now if $x=(z,f(z)) \in N$, then 
 $$
 T_x N = \{ (v, D(f)_zv)\in T_x\Cb^{m+n}: v \in T_z \Cb^n\}.
 $$
Further, by hypothesis, 
$$
(iv, iD(f)_zv)=i \cdot (v,D(f)_zv)   \in T_x N
$$
for all $v \in T_z \Cb^n$. So we must have $iD(f)_z v = D(f)_z iv$ for all $z \in U$ and $v \in T_z \Cb^n$, i.e. $f$ is holomorphic. 
\end{proof}

We may now assemble the above results to prove Lemma \ref{differentiable general}. 

\begin{proof}[Proof of Lemma \ref{differentiable general}]
By Proposition \ref{prop: wedge}, it suffices to prove the lemma for the case when $k=1$. Suppose $\Gamma\subset\PSL(d,\mathbb C)$ is a $(d-1)$-hyperconvex subgroup and $\partial\Gamma$ is a two-sphere. By Corollary~\ref{differentiable k=1} and Lemma \ref{complex submanifold}, $\Lambda_1(\Gamma)$ is a complex submanifold of $\Pb(\Cb^d)$. 

Let $ \Gr_1( T \Pb(\Cb^d))$ denote the 1-Grassmanian bundle  of the tangent bundle of $\Pb(\Cb^d)$, i.e. the fiber above a point $x \in \Pb(\Cb^d)$ consists of all complex lines in $T_x \Pb(\Cb^d)$. Then there is a natural biholomorphism $F : \Fc_{1,2}(\Cb^d) \rightarrow \Gr_1( T \Pb(\Cb^d))$ given by 
$$
F(x) = (x^1, T_{x^1} \Pb(x^2)). 
$$
Since  $\Lambda_1(\Gamma)$ is a complex manifold, 
$$
N : = \{ (x, T_x \Lambda_1(\Gamma) ) : x \in \Lambda_1(\Gamma)\}
$$
is a complex submanifold of  $\Gr_1( T \Pb(\Cb^d))$. By Equation~\eqref{eqn:the tangent spaces}, 
$$
F^{-1}(N) = \Lambda_{1,2}(\Gamma)
$$
and so $\Lambda_{1,2}(\Gamma)$ is a complex submanifold. 
\end{proof}

\subsection*{Sketch of a more dynamical proof} We will now give a sketch of our original, more dynamical proof prior to the simplification suggested by Andres Sambarino. The argument requires that the Zariski closure of $\Gamma$ is simple. 

For any $z\in\Lambda_{\theta_k}(\Gamma)$, Farre--Pozzetti--Viaggi \cite{FPV} defined a homeomorphism
\[\Lambda_k(\Gamma)\to\mathbb P\left(z^{d-k+1}/z^{d-k-1}\right)\cong \mathbb{CP}^1\]
called the \emph{tangent projection at $z$}, and verified that it is bilipschitz on all compact subsets of $\Lambda_k(\Gamma)-\{z^k\}$. Via this homeomorphism, the $\Gamma$-action on $\Lambda_k(\Gamma)$ induces a $\Gamma$-action 
\[\rho_z:\Gamma\to{\rm Homeo}(\mathbb{CP}^1)\]
on $\mathbb{CP}^1$, which they then prove is uniformly quasiconformal. A theorem of Sullivan \cite{sullivan-ergodic} then implies that there exists a a uniform lattice $\hat\Gamma\subset\mathsf{PSL}(2,\mathbb C)$ and a quasi-conformal conjugacy $\phi:\rho_z(\Gamma)\to\hat\Gamma$.

By Theorem \ref{prop: Ahlfors regularity}, the $\alpha_k$-Patterson--Sullivan measure for $\Gamma$ on $\Lambda_k(\Gamma)$ is Ahlfors regular and the $\alpha_1$-Patterson--Sullivan measure of $\hat\Gamma$ on $\mathbb{CP}^1$ is mutually absolutely continuous with the Lebesgue measure. Since quasi-conformal maps preserve the Lebesgue measure class, we can then conclude that their pull-backs to $\partial\Gamma$ are mutually absolutely continuous.
Then a result of Blayac--Canary--Zhu--Zimmer \cite[Thm.\ 13.2]{BCZZ} (see also Sambarino \cite[Cor.\ 3.3.4]{sambarino-dichotomy}) implies that
\begin{equation}
\label{lengths equal}
\alpha_k(\nu(\gamma))=\alpha_1(\nu(\phi(\rho_z(\gamma))))
\end{equation}
for all $\gamma\in\Gamma$, where $\nu$ on the left denotes the Jordan projection in $\PSL(d,\Cb)$, while $\nu$ on the right denotes the Jordan projection for $\PSL(2,\Cb)$. 

We then consider the identity component $\mathsf{G}^0$ of the real Zariski closure $\mathsf{G}$ of $\Gamma$. \textbf{Assuming} that $\mathsf{G}^0$ is simple, we may use an argument of Dal'bo--Kim \cite{DK1} to show that $\mathsf{G}^0$ is isomorphic to $\mathsf{PSL}(2,\mathbb C)$. Once we have
done so, we can complete the proof of Theorem \ref{thm:entropy rigidity strongly irred} just as above. 

Finally, we briefly sketch the argument that $\mathsf{G}^0$ is isomorphic to $\mathsf{PSL}(2,\mathbb C)$. Consider the subgroup
$$H_0:=\{(\gamma,\phi(\rho_x(\gamma)): \gamma\in\Gamma\cap\Gsf^0\}\subset \mathsf{G}^0\times \mathsf{PSL}(2,\mathbb C),$$
and let $Z$ be its real Zariski closure in $\mathsf{G}^0\times \mathsf{PSL}(2,\mathbb C).$
Equation (\ref{lengths equal}) and a result of Benoist \cite{benoist} implies that $H_0$ is not Zariski dense, so
$Z\ne\mathsf{G}^0\times \mathsf{PSL}(2,\mathbb C).$ Consider the projection maps
$\pi_1: Z\to \mathsf{G}^0$ and $\pi_2: Z\to \mathsf{PSL}(2,\mathbb C)$, each of which is surjective. The kernel of $\pi_1$ is  normal in
$\{1\}\times\mathsf{PSL}(2,\mathbb C)$ and cannot be all of $\{1\}\times\mathsf{PSL}(2,\mathbb C)$. Since $\mathsf{PSL}(2,\mathbb C)$ is simple and has trivial center, $\pi_1$ is an isomorphism. Similarly, since $\mathsf{G}^0$ is simple and has trivial center, $\pi_2$ is an isomorphism.
Therefore, $\pi_1\circ\pi_2^{-1}$ is an isomorphism from $\mathsf{G}^0$ is isomorphic to $\mathsf{PSL}(2,\mathbb C)$.

%% file: appendix.tex

\appendix 

\section{Reducing to the projective case}

In this appendix we give the proof of Proposition \ref{prop: wedge}.

For the rest of the section fix $k \in \{1,\dots, d-1\}$ and, as usual, let $(e_1,\dots,e_d)$ denote the standard basis on $\Kb^d$. We identify 
$$
\Kb^{\bar{d}} = \wedge^k \Kb^d, \quad \bar{d} := {d\choose k}=\dim\wedge^k\Kb^d,
$$
using the basis 
\[(e_{i_1}\wedge\dots\wedge e_{i_k})_{1\le i_1<\dots<i_k\le d}\]
where the order is in lexicographic order of the subscripts. This gives an identification $\PSL(\wedge^k\Kb^d)=\PSL(\bar{d},\Kb)$. There is a natural representation $\tilde \eta:\SL(d,\Kb)\to\SL(\wedge^k\Kb^d)=\SL(\bar{d},\Kb)$ defined by 
\[
\tilde \eta(g)(v_1\wedge \dots\wedge v_k)=g(v_1)\wedge\dots\wedge g(v_k),
\] 
which descends to a representation 
\[\eta:\PSL(d,\Kb)\to\PSL(\wedge^k\Kb^d)=\PSL(\bar{d},\Kb).\]

Set $\bar\Delta:=\{1,\dots,\bar d-1\}$ and $\bar\theta_1:=\{1,2,\bar d-2,\bar d-1\}\subset\bar\Delta$. For any non-empty subset $\theta\subset\Delta$ (respectively, $\bar\theta\subset\bar\Delta$), let $\pi_\theta:\Fc_\Delta(\Kb^d)\to\Fc_\theta(\Kb^d)$ (respectively, $\bar\pi_{\bar \theta}:\Fc_{\bar \Delta}(\Kb^{\bar d})\to\Fc_{\bar \theta}(\Kb^{\bar d})$) denote the forgetful projection. We record the following well-known properties of $\eta$ that we will use.
\begin{enumerate}
\item[(I)] The representation $\eta$ sends $\Psf_\Delta\subset\PSL(d,\Kb)$ to $\Psf_{\bar\Delta}\subset\PSL(\bar d,\Kb)$, 
so we may define an $\eta$-equivariant, smooth (holomorphic when $\Kb=\Cb$) embedding
$$E:\Fc_\Delta(\Kb^d)\to\Fc_{\bar \Delta}(\Kb^{\bar d})\quad\text{given by}\quad
E(g\Psf_\Delta)=\eta(g)\Psf_{\bar \Delta}.$$
\item[(II)] For any $x\in\Fc_{\Delta}(\Kb^d)$, we have
\begin{align*}
E(x)^1&=\Span(\{v_1\wedge\dots\wedge v_k:\Span(v_1,\dots,v_k)=x^k\}),\\
E(x)^2&=\Span(\{v_1\wedge\dots\wedge v_k:x^{k-1}\subset\Span(v_1,\dots,v_k)\subset x^{k+1}\}),\\
E(x)^{\bar{d}-2}&=\Span(\{v_1\wedge \dots\wedge v_k:\Span(v_1,\dots,v_k)\cap x^{d-k-1}\neq 0\}\\
&\hspace{1cm}\cup\{v_1\wedge \dots\wedge v_k:\dim(\Span(v_1,\dots,v_k)\cap x^{d-k+1})\ge 2\}),\\
E(x)^{\bar{d}-1}&=\Span(\{v_1\wedge\dots\wedge v_k:\Span(v_1,\dots,v_k)\cap x^{d-k}\neq 0\}).
\end{align*}
In particular, $E(x)^1 = E_k(x^k)$ and 
we have well-defined, $\tau$-equivariant, smooth embeddings
\[E_{\theta_k}:\Fc_{\theta_k}(\Kb^d)\to\Fc_{\bar\theta_1}(\Kb^{\bar d}).\]
for any $x \in \Fc_\Delta(\Kb^d)$.
\item[(III)] The embedding $E$ preserves transversality of pairs of flags, so the same is true for $E_{\theta_k}$.
\item[(IV)] The representation $\tilde \eta$ sends $\mathsf{SU}(d,\Kb)$ into $\mathsf{SU}(\bar d,\Kb)$. 
\item[(V)] If $a={\rm Diag}(a_1,\dots,a_d)\in\PSL(d,\Kb)$ is real valued, then $\tilde \eta(a)$ is a diagonal, real valued matrix whose diagonal entries are all possible products of $k$ of the $d$ diagonal entries of $a$, and whose $(1,1)$ and $(2,2)$ entries are 
\[\prod_{i=1}^ka_i\quad\text{and}\quad a_{k+1}\prod_{i=1}^{k-1}a_i\]
respectively. Furthermore, if $a_1\ge \dots \ge a_d$, then the largest and second largest
entries of $\tilde \eta(a)$ are respectively its $(1,1)$ and $(2,2)$ 
entries, and its third largest entry is 
\[\left\{\begin{array}{ll}
\displaystyle a_{k+1}a_k\prod_{i=1}^{k-2}a_i&\text{if }d=k+1,\\
\displaystyle \max\left\{a_{k+2}\prod_{i=1}^{k-1}a_i,\,a_{k+1}a_k\prod_{i=1}^{k-2}a_i\right\}&\text{if }d\ge k+2.\end{array}\right.\]
\item[(VI)] By (I), (IV), and (V), if $g=man$ is the Iwasawa decomposition of $g\in\PSL(d,\Kb)$, then $\eta(g)=\eta(m)\eta(a)\eta(n)$ is the Iwasawa decomposition of $\eta(g)$, and $\alpha_1(\log \tau(a))=\alpha_k(\log a)$.
In particular, if $\bar B:\PSL(\bar d,\Kb)\times \Fc_{\bar \Delta}\to\Rb$ is the Iwasawa cocycle of $\PSL(\bar d,\Kb)$, then for all $g\in\PSL(d,\Kb)$ and $x\in\Fc_\Delta(\Kb^d)$, we have
\[\alpha_1(\bar B(\eta(g),E(x))=\alpha_k(B(g,x)).\]
\item[(VII)] By (IV) and (V), if $g=m_1 e^{\kappa(g)} m_2$ is a singular value decomposition of $g$, then there is a permutation matrix $p\in\mathsf{PGL}(\bar d,\Kb)$ whose $(1,1)$, $(2,2)$, $(\bar d-1,\bar d-1)$, and $(\bar d,\bar d)$ entries are all $1$ such that $\eta(g)=\bar m_1e^{\kappa(\eta(g))}\bar m_2$ is a singular value decomposition of $\eta(g)$, where $\bar m_1:=\eta(m_1)p^{-1}$ and $\bar m_2:=p\eta (m_2)$. From this, we may deduce that $\alpha_1(\kappa(\eta(g)))=\alpha_k(\kappa(g))$,
\[\alpha_2(\kappa(\eta(g)))=\left\{\begin{array}{ll}
\alpha_{k-1}(\kappa(g))&\text{if }d=k+1,\\
\max\{\alpha_{k-1}(\kappa(g)),\alpha_{k+1}(\kappa(g))\}&\text{if }d\ge k+2,
\end{array}\right.\]
and $E_{\theta_k}(\Usf_{\theta_k}(g))=\Usf_{\bar\theta_1}(\eta(g))$. 
\end{enumerate}

\begin{proof}[Proof of Proposition~\ref{prop: wedge}]
Notice that $E(x)^1$ and $E(x)^2$ only depend on $x^{k-1}, x^k, x^{k+1}$. Hence there exists a unique map $E_k: \Fc_{k-1,k,k+1}(\Kb^d) \rightarrow \Fc_{1,2}(\Kb^{\bar d})$ satisfying 
$$
E_k(x^{k-1}, x^k, x^{k+1}) = ( E(x)^1, E(x)^2)
$$
for all $x \in \Fc_{\theta_k}(\Kb^d)$. Moreover, the formula for $E(x)^1$ and $E(x)^2$ imply that $E_k$ is a $\eta_k$-equivariant smooth embedding, which is holomorphic when $\Kb=\Cb$. 

The formula for $E(x)^1$ implies that $E_k$ descends to a smooth embedding $E_k^1 : \Gr_k(\Kb^d) \rightarrow \Pb(\Rb^{d_k})$. Hence (3) is true.

Next we prove (1). By (VII), the fact that $\Gamma$ is $\Psf_{\theta_k}$-divergent implies that $\eta(\Gamma)$ is $\Psf_{\bar\theta_1}$-divergent, 
and $\Lambda_{\bar\theta_1}(\eta(\Gamma))=E_{\theta_k}(\Lambda_{\theta_k}(\Gamma))$. In particular, 
$\Lambda_1(\eta(\Gamma))=E_k(\Lambda_k(\Gamma))$. Since $\Gamma$ is $\Psf_{\theta_k}$-transverse, (III) implies that $\eta(\Gamma)$ is 
$\Psf_{\bar\theta_1}$-transverse. Since $E_{\theta_k}$ is an $\eta$-equivariant embedding, the fact that every point in 
$\Lambda_{\theta_k}(\Gamma)$ is a conical point for the $\Gamma$-action implies that every point in 
$\Lambda_{\bar\theta_1}(\eta(\Gamma))$ is a conical point for the $\tau(\Gamma)$-action. Hence, $\eta(\Gamma)$ is $\Psf_{\theta_k}$-Anosov.

Let $x,y,z\in \Lambda_{\theta_k}(\Gamma)$ be any mutually distinct triples of flags. Since $x$ and $z$ are transverse, we may choose a basis $(v_1,\dots,v_d)$ of $\Kb^d$ such that 
\[x^{k-1}=\Span(v_1,\dots,v_{k-1}),\quad x^k\cap z^{d-k+1}=\Span(v_k),\]
 \[x^{k+1}\cap z^{d-k}=\Span(v_{k+1}),\quad\text{and}\quad z^{d-k-1}=\Span(v_{k+2},\dots,v_d).\] 
 For any tuple ${\bf j}=(j_1,\dots,j_k)$ of integers such that $1\le j_1<\dots<j_k\le d$, set $v_{\bf j}:=v_{j_1}\wedge\dots\wedge v_{j_k}$. 
 Then observe from the description of $E_{\theta_k}$ given in (II) that 
 \[E_{\theta_k}(x)^1=\Span(v_{(1,\dots,k)})\quad\text{and}\quad E_{\theta_k}(z)^{\bar d-2}=\Span(\{v_{\bf j}:{\bf j}\neq (1,\dots,k),(1,\dots,k-1,k+1)\}).\]
 Also, since $\Gamma$ is $(d-k)$-hyperconvex, we have
$$\Big((x^k\cap z^{d-k+1})+z^{d-k-1}\Big)\cap \Big( (y^k\cap z^{d-k+1})+z^{d-k-1}\Big)=z^{d-k-1}.$$
Since $v_k \in (x^k\cap z^{d-k+1})+z^{d-k-1}$ and $v_k \notin z^{d-k-1}$, we must have 
$$
v_k \notin (y^k\cap z^{d-k+1})+z^{d-k-1}.
$$
Hence, there is some $a_k,a_{k+2},\dots,a_d\in\Rb$ such that 
\[y^k\cap z^{d-k+1}=\Span\left(v_{k+1}+a_kv_k+\sum_{j=k+2}^da_jv_j\right).\]
Since $y^{k-1}+z^{d-k+1}=\Kb^d$, it now follows that 
 \[E_{\theta_k}(y)^1=\Span\left(v_{(1,\dots,k-1,k+1)}+\sum_{{\bf j}\neq (1,\dots,k-1,k+1)}b_{\bf j}v_{\bf j}\right)\]
 for some $\{b_{\bf j}\}\subset\Kb$. Thus, $E_{\theta_k}(x)^1+E_{\theta_k}(y)^1+E_{\theta_k}(z)^{\bar d-2}=\Kb^{\bar d}$. 
 Since $x,y,z$ was arbitrary, it follows that $\eta(\Gamma)$ is $(\bar{d}-1)$-hyperconvex. This finishes the proof of (1).
 
Statements (2) and (4) are immediate consequences of (VII) and (VI) respectively.
\end{proof}